\documentclass[11pt]{amsart}
\usepackage{geometry}                
\usepackage[parfill]{parskip}    
\usepackage{graphicx}
\usepackage{amssymb}
\usepackage{epstopdf}

\title{Doing and Showing}
\author{Andrei Rodin}

\begin{document}
\maketitle

\begin{abstract}
The persisting gap between the formal and the informal
mathematics is due to an inadequate notion of mathematical
theory behind the current formalization techniques. I mean the
(informal) notion of axiomatic theory according to which a
mathematical theory consists of a set of axioms and further
theorems deduced from these axioms according to certain rules of
logical inference. Thus the usual notion of axiomatic
method is inadequate and needs a replacement.
\end{abstract}

\renewcommand{\contentsname}{}
\tableofcontents
\section{Introduction}
\paragraph{An  \emph{axiomatic theory} consists of a distinguished set of propositions called  \emph{axioms} and another set of propositions called  \emph{theorems}, which are deduced from the axioms according to certain  \emph{rules of inference}. Such rules are supposed to be \emph{truth-preserving} in the following sense: as far as the axioms are true the theorems derived from these axioms are also true. Further, these rules are supposed to be not specific for any given theory: one assumes that the same set of rules of inference applies in all axiomatic theories (logical monism) or at least that any complete set of such rules applies in  some large class of theories (logical pluralism). These basic features of the rules of inference can be briefly expressed by saying that those rules are the rules of \emph{logical} inference. A logical inference of a given theorem from the axioms (possibly through a number of intermediate propositions) is called a \emph{proof} of this theorem.
}

\paragraph{
Here is how this core notion of axiomatic theory is described by Hilbert in his famous address ``Axiomatic Thought'' delivered before the Swiss Mathematical Society in Zurich in 1917:
}

\begin{quote}
If we consider a particular theory more closely, we always
see that a few distinguished propositions of the field of
knowledge underlie the construction of the framework of
concepts, and these propositions then suffice by
themselves for the construction, in accordance with
logical principles, of the entire framework. ...  These fundamental propositions can be regarded ... as
the \underline{axioms} of the individual fields of knowledge : the
progressive development of the individual field of
knowledge then lies solely in the further logical
construction of the already mentioned framework of
concepts. This standpoint is especially predominant in
pure mathematics. ... [A]nything at all that can be the object of scientific
thought becomes dependent on the axiomatic method,
an thereby indirectly on mathematics. (\cite{Hilbert:1918})
\end{quote}

\paragraph{In a different paper \cite{Hilbert:1922} Hilbert goes even further and claims that: 
}

\begin{quote}
The axiomatic method is and remains the indispensable
tool, appropriate to our minds, for all exact research in
any field whatsoever : it is logically incontestable and at
the same time fruitful. ... To proceed axiomatically means
in this sense nothing else than to think with
consciousness.
\end{quote}

\paragraph{The above quotes do not fully describe what Hilbert has in mind when he talks about the axiomatic method; this general description leaves a wide room for further specifications and interpretations, many of which turn to be mutually incompatible. Discussions over such further details of the axiomatic method of theory-building have been playing an important role in the philosophy of mathematics since the beginning of the 20th century. One such continuing discussion concerns the epistemic status of axioms. According to the traditional view dating back to Aristotle and more recently defended by Frege, axioms are self-evident truths, which can and should be used for proving some further propositions (theorems), which by themselves are not evident. According to a novel viewpoint defended by Hilbert, axioms are formal expressions that given different truth-values through their different \emph{interpretations}  (I provide some further explanations in Section 6 below). Such different views on axioms lead Frege and Hilbert to very different notions of axiomatic theory  and different understanding of the axiomatic method. However important this and other differences in the understanding of the core axiomatic method might be they are not directly relevant to my argument given below in this paper. Instead of discussing various specifications of the axiomatic method I would like to put here under a critical examination the core notion of axiomatic method itself (which is assumed, in particular, by both Hilbert and Frege in their influential debate \cite{Frege:1971}).  
}   

\paragraph{
The paper is organized as follows. First, I elaborate  in some detail on the First Book of Euclid's ``Elements'' and show that Euclid's theory of geometry is \underline{not} axiomatic in the modern sense but is construed differently. Second, I provide some evidences showing that the usual commonly accepted notion of axiomatic theory equally fails to account for today's mathematical theories. I also provide some polemical arguments against the popular view according to which a good mathematical theory \emph{must} be axiomatic (in the usual general sense) and point to an alternative method of theory-building. Since my critique of the core axiomatic method is \emph{constructive} in its character I briefly observe known constructive approaches in the foundations of mathematics and describe the place of my proposal in this context. The main difference of my and earlier constructive proposals for foundations of mathematics appears to be the following: while earlier proposals deal with the issue of admissibility of some particular mathematical principles and like choice and some putative mathematical objects like infinite sets my proposal concerns the very method of theory-building. As a consequence, my proposal unlike earlier constructive proposals puts no restriction on the existing mathematical practice but rather suggests an alternative method of organizing this practice into a systematic theoretical form. In the concluding section of the paper I argue that the constructive mathematics (in the specific sense of the term specified in an earlier section of this paper) better serves needs of mathematically-laden empirical sciences than the formalized mathematics.          
}

\part{Euclid's Way of Building Mathematical Theories}
\section{Demonstration and ``Monstration''} 

\paragraph{
All Propositions of Euclid's \emph{Elements} \cite{Euclides:1883-1886} (with few easily understandable exceptions) are designed according to the same scheme described by Proclus in his \emph{Commentary} \cite{Proclus:1970} as follows: 
}
\begin{quote}
Every Problem and every Theorem that is furnished with all its parts should contain the
following elements: an \emph{enunciation}, an \emph{exposition}, a \emph{specification}, a
\emph{construction}, a \emph{proof}, and a \emph{conclusion}. Of these \emph{enunciation} states what is given and
what is being sought from it, a perfect \emph{enunciation} consists of both these parts. The \emph{exposition}
takes separately what is given and prepares it in advance for use in the investigation. The
\emph{specification} takes separately the thing that is sought and makes clear precisely what it is. The
\emph{construction} adds what is lacking in the given for finding what is sought. The \emph{proof} draws the
proposed inference by reasoning scientifically from the propositions that have been admitted.
The \emph{conclusion} reverts to the \emph{enunciation}, confirming what has been proved. (\cite{Proclus:1970}, p.203, italic is mine)
\end{quote}

\paragraph{
It is appropriate to notice here that the term ``proposition'', which is traditionally used in translations  as a common name of Euclid's problems and theorems is not found in the original text of the \emph{Elements}: Euclid numerates these things throughout each Book without naming them by any common name. (The reader will shortly see that this detail is relevant.) The difference between problems and theorems will be explained in the Section 4 below. Let me now show how this Proclus' scheme applies to Proposition 5 of the First Book (Theorem 1.5), which is a well-known theorem about angles of the isosceles triangle. I quote Euclid's \emph{Elements} by the recent English translation by Richard Fitzpatrick \cite{Euclid:2011}. References in square brackets are added by the translator; some of them will be discussed later on. Words in round brackets are added by the translator for stylistic reason. Words in angle brackets are borrowed from the above Proclus' quote. Throughout this paper I write these words in italic when I use them in Proclus' specific sense. 
}

\paragraph{[\emph{enunciation:}]}
\begin{quote}
For isosceles triangles, the angles at the base are equal to one another, and if the equal straight lines are produced then the angles under the base will be equal to one another.
\end{quote}

\paragraph{[\emph{exposition}]:}

\begin{quote}
Let $ABC$ be an isosceles triangle having the side $AB$ equal to the side $AC$; and let the straight lines $BD$ and $CE$ have been produced further in a straight line with $AB$ and $AC$ (respectively). [Post. 2].
\end{quote}

\includegraphics [scale=0.5]{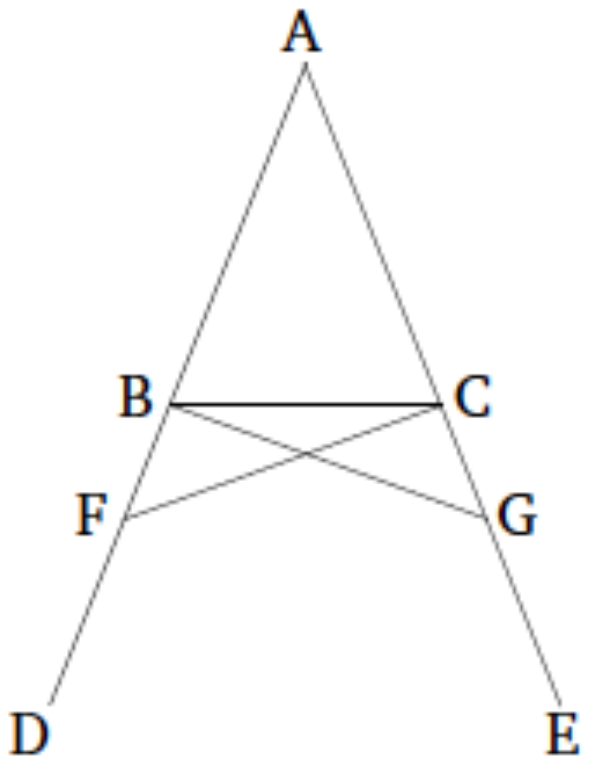}

\paragraph{[\emph{specification}:]}

\begin{quote}
I say that the angle $ABC$ is equal to $ACB$, and (angle) $CBD$ to $BCE$.
\end{quote}

\paragraph{[\emph{construction}:]}

\begin{quote}
For let a point $F$ be taken somewhere on $BD$, and let $AG$  have been cut off from the greater  $AE$,  equal to the lesser $AF$ [Prop. 1.3]. Also, let the straight lines $FC$, $GB$ have been joined. [Post. 1] 
\end{quote}

\paragraph{[\emph{proof}:]}

\begin{quote}
In fact, since $AF$ is equal to $AG$, and $AB$ to $AC$,
the two (straight lines) $FA$, $AC$ are equal to the two (straight lines) $GA$, $AB$, respectively. They also encompass a common angle $FAG$.
Thus, the base $FC$ is equal to the base $GB$, and the triangle $AFC$ will be equal to the triangle $AGB$,
and the remaining angles subtended by the equal sides will be equal  to the corresponding remaining angles [Prop. 1.4]. (That is) $ACF$ to $ABG$, and $AFC$ to $AGB$.
And since the whole of $AF$ is equal to the whole of $AG$,
within which $AB$ is equal to $AC$, the remainder $BF$ is thus equal to the remainder $CG$ [Ax.3].
But $FC$ was also shown (to be) equal to $GB$.
So the two (straight lines) $BF$, $FC$ are equal to the two (straight lines) $CG$, $GB$ respectively, and the angle $BFC$ (is) equal to the angle $CGB$,
while the base $BC$ is common to them. Thus
the triangle $BFC$ will be equal to the triangle $CGB$,
and the remaining angles subtended by the equal sides will be equal to the corresponding remaining angles [Prop. 1.4]. Thus $FBC$ is equal to $GCB$, and $BCF$ to $CBG$. Therefore, since the whole angle $ABG$ was shown (to be) equal to the whole angle $ACF$, within which 
$CBG$ is equal to $BCF$,
the remainder $ABC$ is thus equal to the remainder $ACB$ [Ax. 3].
And they are at the base of triangle $ABC$.
And $FBC$ was also shown (to be) equal to $GCB$.
And they are under the base.\end{quote}

\paragraph{[\emph{conclusion}:]}

\begin{quote}
Thus, for isosceles triangles, the angles at the base are
equal to one another, and if the equal sides are produced
then the angles under the base will be equal to one another.
(Which is) the very thing it was required to show.
\end{quote}

\paragraph{An obvious difference between Proclus' analysis of the above theorem and its usual modern analysis is the following. For a modern reader the proof of this theorem begins with Proclus' \emph{exposition} and includes Proclus' \emph{specification}, \emph{construction} and \emph{proof}. Thus for Proclus the \emph{proof} is only a part of what we call today the proof of this theorem. Also notice that Euclid's theorems conclude with the words ``which ... was required to \emph{show}'' (as correctly translates Fitzpatrick) but not with the words ``what it was required to \emph{prove}'' (as inaccurately translates Heath \cite{Heath:1926}). The standard Latin translation of this Euclid's formula as \emph{quod erat demonstrandum} is also inaccurate. These inaccurate translations conflate two different Greek verbs: ``apodeiknumi'' (English ``to prove'', Latin ``demonstrare'') and ``deiknumi'' (English ``to show'', Latin ``monstrare''). The difference between the two verbs can be clearly seen in the two Aristotle's \emph{Analytics}: Aristotle uses the verb ``apodeiknumi'' and the derived noun ``apodeixis'' (proof) as technical terms in his syllogistic logic, and he uses the verb ``deiknumi'' in a broader and more informal sense when he discusses epistemological issues (mostly in the \emph{Second Analytics}). Without trying to trace here the history of Greek logical and mathematical terminology and speculate about possible influences of some Greek writers on some other writers I would like to stress the remarkable fact that Aristotle's use of verbs ``deiknumi'' and ``apodeiknumi'' agrees with Euclid's and Proclus'. In my view this fact alone is  sufficient for taking  seriously the difference between the two verbs and distinguish between \emph{proof} and ``showing'' (or otherwise between \emph{demonstration} and \emph{monstration}).   
}

\footnote{As far as mutual influences are concerned two things are certain: (i) Proclus read Aristotle and (ii) Aristotle had at least a basic knowledge of the mathematical tradition, on which Euclid later elaborated in his \emph{Elements} (as Aristotle's mathematical examples clearly show \cite{Heath:1949}).  It remains unclear whether  Aristotle's work could influence Euclid. In my view this is unlikely. However Aristotle's logic certainly played an important role in later interpretations and revisions of Euclid's \emph{Elements}. I leave this interesting issue outside of the scope of this paper.}

\paragraph{
One may think that the difference between the current meaning of the word ``proof'' in today's mathematics and logic and the meaning of Proclus'  \emph{proof} (Greek ``apodeixis'') is a merely terminological issue, which is due to difficulties of translation from Greek to English. I shall try now to show that this terminological difference points on a deeper problem, which is not merely linguistic. In today's logic the word ``proof'' stands for a logical inference of certain conclusion from some given premises. In fact this is what by and large was meant by proof also by Aristotle and Proclus. Indeed, looking at the \emph{proof} (in Proclus' sense) of Euclids Theorem 1.5 we see that it also qualifies as a proof in the modern sense: we have here a number of premises (which I shall specify in the next Section) and certain conclusions derived from those premises. It is irrelevant now whether or not this particular inference is valid according to today's logical standards; what I want to stress here is only the general setting that involves some premises, an inference (probably invalid) and some conclusions. This core meaning of the word ``proof'' (Greek ``apodeixis'') hardly changed since Proclus' times.}

\paragraph{So we get a problem, which is clearly not only terminological: Is it indeed justified to describe the \emph{exposition}, the \emph{specification} and the \emph{construction} as elements of the proof or one should rather follow Proclus and consider these things as independent constituents of a mathematical theorem? }

\paragraph{The question of \emph{logical significance} of the \emph{exposition}, the \emph{specification} and the \emph{construction} in Euclid's geometry has been discussed in the literature; in what follows I shall briefly describe some tentative answers to it. However before doing this I would like to stress that this question may be ill-posed to begin with. As far as one assumes, first, that the theory of Euclid's  \emph{Elements} is  (by and large) sound and, second, that any sound mathematical theory is an axiomatic theory in the modern sense, then, in order to make these two assumptions mutually compatible, one has to describe the \emph{exposition}, the \emph{specification} and the \emph{construction} of each Euclid's theorem as parts of the proof of this theorem and specify their logical role and their logical status. I shall not challenge the usual assumption according to which Euclid's mathematics is by and large sound. (I say ``by and large'' in order to leave some room for possible revisions and corrections of Euclid's arguments and thus avoid controversies about the question whether a given interpretation of Euclid is authentic or not. Although I pay more attention to textual details than it is usually done in modern logical reconstructions of Euclid's reasoning, I am not going to criticize these reconstructions by pointing to their anachronistic character.) However I shall challenge the other assumption according to which any sound mathematical theory is an axiomatic theory in the modern sense. Since I do not take this latter assumption for granted I do not assume from the outset that the problematic elements of Euclid's reasoning  (the \emph{exposition}, the \emph{specification} and the \emph{construction}) play some \emph{logical} role, which only needs to be made explicit and appropriately understood. In what follows I try to describe how these elements work without making about them any additional assumptions and only then decide whether the role of these elements qualifies as logical or not.   
}

\section{Are Euclid's Proofs Logical?}
Let's look at Euclid's Theorem 1.5 more attentively. I begin its analysis with its \emph{proof}. Among the premisses of this \emph{proof}, one may easily identify Axiom (Common Notion) 3 according to which

\begin{quote}
(Ax.3): If equal things are subtracted from equal things then the remainders are equal
\end{quote}

and the preceding Theorem 1.4 according to which

\begin{quote}
(Prop.1.4): If two triangles have two corresponding sides equal, and have the angles enclosed by the equal sides equal, then they will also have equal bases, and the two triangles will be equal, and the remaining angles subtended by the equal sides will be equal to the corresponding remaining angles.
\end{quote}

 \paragraph{
 I shall not comment on the role Theorem 1.4 in this \emph{proof} (which seems to be clear) but say few things about the role of the Axiom 3.
}
Here is how exactly the Axiom (Common Notion) 3 is used in the above Euclid's \emph{proof}. First, \emph{by construction} we have

\begin{quote}
\textbf{Con1}: $BF \equiv AF - AB$ and    
\textbf{Con2}: $CG \equiv AG - AC$
\end {quote}

which is tantamount to saying that point $B$ lays between points $A$, $F$ and point $C$ lays between points $A$, $G$). Second, \emph{by hypothesis} we have

\begin{quote}
\textbf{Hyp}: $AB = AC$ 
\end{quote}

and once again \emph{by construction}

\begin{quote}
\textbf{Con3}: $AF = AG$ 
\end{quote}

\paragraph{
Now we see that we have got the situation described in Ax.3: equal things are subtracted from equal things. Using this Axiom we conclude that $BF = CG$.
}
 
\paragraph{Notice that Ax.3 applies to all ``things'' (mathematical objects), for which the relation of \emph{equality} and the operation of \emph{subtraction} make sense. In Euclid's mathematics this relation and this operation apply not only to straight segments and numbers but also to geometrical  objects of various sorts including \emph{figures}, angles and solids. Since Euclid's equality is not interchangeable with identity I use for the two relations two different symbols: namely I use the usual symbol for Euclid's equality (even if this equality is not quite usual), and use symbol $\equiv$ for identity. My use of symbols $+$ and $-$ is self-explanatory.} 

\footnote{ The \emph{difference} $A - B$ of two figures $A$, $B$ is a figure obtained through ``cutting'' $B$ out of $A$; the  \emph{sum}  $A + B$  is the result of \emph{concatenation} of  $A$ and $B$. These operations are not defined up to \emph{congruence} of figures (for there are, generally speaking, many possible ways, in which one may cut out one figure from another) but, according to Euclid's Axioms, these operations are defined up to Euclid's \emph{equality}. This shows that Euclid's \emph{equality} is weaker than \emph{congruence}: according to Axiom 4 congruent objects are equal but, generally, the converse does not hold. In the case of (plane) figures Euclid's equality is equivalent to the equality (in the modern sense) of their air.     
}

The other four Euclid's Axioms (not to be confused with Postulates!) have  the same character.
This makes Euclid's Axioms in general, and Ax.3 in particular, very unlike premises like \textbf{Con1-3} and \textbf{Hyp}, so one may wonder whether the very idea of treating these things on equal  footing (as different \emph{premises} of the same inference) makes sense. More precisely we have here the following choice.
 One option is to interpret Ax.3 as the following implication:

\begin{quote}
$\{(a \equiv b - c) \& (d \equiv e - f ) \& (b = d ) \& (c =f)\} \rightarrow (a = b)$
\end{quote} 

\paragraph{and then use it along with \textbf{Con1-3} and \textbf{Hyp} for getting the desired conclusion through \emph{modus ponens} and other appropriate rules. This standard analysis involves a fundamental distinction between premises and conclusion, on the one hand, and rules of inference, on the other hand. It assumes that in spite of the fact that Euclid (as most of other mathematicians of all times) remains silent about logic, his reasoning nevertheless follows some implicit logical rules. The purpose of logical analysis  in this case is to make this ``underlying logic`` (as some philosophers like to call it) explicit.}
  
\paragraph{The other option that I have in mind is to interpret Ax.3 itself as a rule rather than as a premiss. Following this rule, which can be pictures as follows:}

\begin{quote}
$(a \equiv b - c), (d \equiv e - f ), (b = d ), (c =f)$\\
------------------------------------------------------------------------ (Ax.3)\\
(a = b)
\end{quote}

\paragraph{one derives from \textbf{Con1-3} and \textbf{Hyp}  the desired conclusion. So interpreted Ax.3 hardly qualifies as a  \emph{logical} rule because it applies only to propositions of a particular sort (namely, of the form $X = Y$ where $X, Y$ are some \emph{mathematical} objects of appropriate types). This Axiom cannot help one to prove that Socrates is mortal. Nevertheless the domain of application of this rule is quite vast and covers the whole of Euclid's mathematics. An important advantage of this analysis is that it doesn't require one to make any assumption about hidden features of Euclid's thinking: unlike the distinction between logical rules and instances of applications of these rules the distinction between axioms and premises like \textbf{Con1-3} and \textbf{Hyp} is explicit in Euclid's \emph{Elements}.}

 \paragraph{There is also a historical reason to prefer the latter reading of Euclid's Common Notions. Aristotle uses the word ``axiom'' interchangeably with the expressions ``common notions'', ``common opinions'' or simply ``commons'' for what we call today logical laws or logical principles but not for what we call today axioms. Moreover in this context he systematically draws an analogy between mathematical common notions and his proposed logical principles (laws of logic). This among other things provides an important historical justification for  calling Euclid's Common Notions by the name of Axioms. It is obvious that mathematics in general and mathematical common notions (axioms) in particular serve for Aristotle as an important source for developing the very idea of logic. Roughly speaking Aristotle's thinking, as I understand it, is this: behind the basic principles of mathematical reasoning spelled out through mathematical common notions (axioms) there are other yet more general principles relevant to reasoning about all sorts of beings and not only about mathematical objects. The fact that Euclid, according to the established chronology, is younger than Aristotle for some 25 years (Euclid's dates unlike Aristotle's are only approximate) shouldn't confuse one. While there is no strong evidence of the influence of Aristotle's work on Euclid, the influence on Aristotle of the same mathematical tradition, on which Euclid elaborated, is clearly documented in Aristotle's writings themselves. In particular, Aristotle quotes Euclid's Ax.3 (which, of course, Aristotle could know from another source) almost verbatim.   
 }
 
 \footnote{Here are some quotes: 
 \begin{quote}
 By first principles of proof [as distinguished from first principles in general] I mean the
common opinions on which all men base their demonstrations, e.g. that one of two
contradictories must be true, that it is impossible for the same thing both be and not to be, and
all other propositions of this kind." (Met. 996b27-32, Heath's translation, corrected)
\end{quote}
Here Aristotle refers to a logical principle as ``common opinion''. In the next quote he compares mathematical and logical axioms:
\begin{quote}
We have now to say whether it is up to the same science or to different sciences to inquire
into what in mathematics is called axioms and into [the general issue of] essence. Clearly the
inquiry into these things is up to the same science, namely, to the science of the philosopher.
For axioms hold of everything that [there] is but not of some particular genus apart from
others. Everyone makes use of them because they concern being qua being, and each genus is.
But men use them just so far as is sufficient for their purpose, that is, within the limits of the
genus relevant to their proofs. Since axioms clearly hold for all things qua being (for being is
what all things share in common) one who studies being qua being also inquires into the
axioms. This is why one who observes things partly [=who inquires into a special domain]
like a geometer or a arithmetician never tries to say whether the axioms are true or false.
(Met. 1005a19-28, my translation)
\end{quote}
Here is the last quote where Aristotle refers to Ax.3 explicitly:
\begin{quote}
Since the mathematician too uses common [axioms] only on the case-by-case basis, it must
be the business of the first philosophy to investigate their fundamentals. For that, when equals
are subtracted from equals, the remainders are equal is common to all quantities, but
mathematics singles out and investigates some portion of its proper matter, as e.g. lines or
angles or numbers, or some other sort of quantity, not however qua being, but as [...]
continuous. (Met. 1061b, my translation)
\end{quote}  
The ``science of philosopher'' otherwise called the ``first philosophy'' is Aristotle's logic, which in his understanding is closely related to (if not indistinguishable from) what we call today ontology. After Alexandrian librarians we call today the relevant collection of Aristotle's texts by the name of \emph{metaphysics} and also use this name for a relevant philosophical discipline. 
}

\paragraph{However important this Aristotle's argument in the history of Western thought may be I see no reason to take it for granted every time when we try today to interpret Euclid's \emph{Elements} or any other old mathematical text. Whatever is one's philosophical stance concerning the place of logical principles in human reasoning one can see what kind of harm can be made if Aristotle's assumption about the primacy of logical and ontological principles is taken straightforwardly and uncritically: one treats Euclid's Axioms on equal footing with premisses like \textbf{Con1-3} and \textbf{Hyp} and so misses the law-like character of the Axioms. Missing this feature doesn't allow one to see the relationships between Greek logic and Greek mathematics, which I just sketched.}

\paragraph{
Having said that I would like to repeat that  Euclid's  \emph{proof} (apodeixis) is the part of Euclid's theorems, which more resembles what we today call proof  (in logic) than other parts Euclid's theorems. For this reason in what follows I shall call inferences in Euclid's \emph{proofs}, which are based on Axioms, \emph{protological} inferences and distinguish them from inferences of another type that I shall call \emph{geometrical} inferences. This analysis is not incompatible with the idea (going back to Aristotle) that behind Euclid's protological and geometrical inferences there are inferences of a more fundamental sort, that can be called \emph{logical} in the proper sense of the word. However I claim that Euclid's text as it stands provides us with no evidence in favor of this strong assumption. One can learn Euclid's mathematics and fully appreciate its rigor without  knowing anything about logic just like Moliere's  M. Jourdain could well express himself long before he learned anything about prose!      
}

\paragraph{
Whether or not the science of logic really helps one to improve on mathematical rigor - or this is rather the mathematical rigor that helps one to do logic rigorously -  is a controversial question that I postpone until the Second Part of this paper.  The purpose of my present reading of Euclid is at the same time more modest and more ambitious than the purpose of logical analysis. It is more modest because this reading doesn't purport to assess Euclid's reasoning from the viewpoint of today's mathematics and logic but aims at reconstructing this reasoning in its  authentic archaic form. It is more ambitious because it doesn't take the today's viewpoint  for granted but aims at reconsidering this viewpoint by bringing it into a historical perspective.    
}

\section{Instantiation and Objectivity}

\paragraph{Let us now see where the premises \textbf{Hyp} and \textbf{Con 1-3} come from. As I have already mentioned they actually come from two different sources: \textbf{Hyp} is assumed  \emph{by hypotheis} while \textbf{Con 1-3} are assumed \emph{by construction}. 
I shall consider here these two cases one after the other.}

\paragraph{The notion of hypothetic reasoning is an important extension of the core notion of axiomatic theory described in Section 1 above; it is well-treated in the literature and I shall not cover it here in full. I shall consider only one particular aspect of hypothetical reasoning as it is present in Euclid. The hypothesis that validates \textbf{Hyp}, informally speaking, amounts to the fact that Theorem 1.5 tells us something about isosceles triangles (rather than about objects of another sort). The corresponding definition (Definition 1.20) tells us that two sides of the isosceles triangle are equal. However to get from here to \textbf{Hyp} one needs yet another step. The \emph{enunciation} of Theorem 1.5 refers to isosceles triangles  \emph{in general}. But \textbf{Hyp} that is involved into the \emph{proof} of this Theorem concerns only \emph{particular} triangle $ABC$. Notice also that the \emph{proof} concludes with the propositions $ABC = ACB$ and $FBC = GCB$ (where $ABC$,  $ACB$, $FBC$ and $GCB$ are angles), which also concern only  \emph{particular} triangle $ABC$. This conclusion differs from the following \emph{conclusion} (of the whole Theorem), which almost verbatim repeats the  \emph{enunciation} and once again refers to isosceles triangles and their angles in general terms.}

\paragraph{The wanted step that allows Euclid to proceed from the \emph{enunciation} to  \textbf{Hyp} is made in the \emph{exposition} of this Theorem, which introduces triangle $ABC$ as an ``arbitrary representative'' of isosceles triangles (in general). In terms of modern logic this step can be described as the \emph{universal instantiation}:\\  
$$\forall x P(x) \Longrightarrow P(a/x)$$\\ 
where $P(a/x)$ is the result of the substitution of individual constant $a$ at the place of all free occurrences of variable $x$ in $P(x)$. The same notion of universal instantiation allows for interpreting Euclid's  \emph{specification} in the obvious way. The reciprocal backward step that allows Euclid to obtain the \emph{conclusion} of the Theorem from the conclusion of the \emph{proof} can be similarly described as the  \emph{universal generalization }:\\  
$$P(a) \Longrightarrow \forall x P(x) $$\\ 
(which is a valid rule only under certain conditions that I skip here).  
}

\paragraph{As long as the \emph{exposition} and the \emph{specification} are interpreted in terms of the universal instantiation these operations are understood as logical inferences and, accordingly, as element of proof in the modern sense of the word. A somewhat different - albeit not wholly incompatible - interpretation of Euclid's  \emph{exposition} and \emph{specification} can be straightforwardly given in terms of Kant's \emph{transcendental aesthetics} and \emph{transcendental logic} developed in his \emph{Critique of Pure Reason} \cite{Kant:1999}. Kant thinks of the traditional geometrical \emph{exposition} not as a logical inference of one proposition from another but as a ``general procedure of the imagination for providing a concept with its image''; a representation of such a general procedure Kant calls a \emph{schema} of the given concept (A140). Thus for Kant any individual mathematical object (like triangle $ABC$)  always comes with a specific \emph{rule} that one follows constructing this object in one's imagination and that provides a link between this object and its corresponding concept (the concept of isosceles triangle in our example). According to Kant the representation of general concepts by imaginary individual objects (which Kant for short also describes as ``construction of concepts'') is the principal distinctive feature of mathematical thinking, which distinguishes it from a philosophical speculation. 
}   

\begin{quote}
``Philosophical cognition is rational cognition from concepts, mathematical cognition is that from the construction of concepts.'' But to construct a concept means to exhibit a priori the intuition corresponding to it. For the construction of a concept, therefore, a non-empirical intuition is required, which consequently, as intuition, is an individual object, but that must nevertheless, as the construction of a concept (of a general representation), express in the representation universal validity for all possible intuitions that belong under the same concept, either through mere imagination, in pure intuition, or on paper, in empirical intuition.... The individual drawn figure is empirical, and nevertheless serves to express the concept without damage to its universality, for in the case of this empirical intuition we have taken account only of the action of constructing the concept, to which many determinations, e.g., those of the magnitude of the sides and the angles, are entirely indifferent, and thus we have abstracted from these differences, which do not alter the concept of the triangle. \\
Philosophical cognition thus considers the particular only in the universal, but mathematical cognition considers the universal in the particular, indeed even in the individual...
 (A713-4/B741-2). 
\end{quote}

\paragraph{Kant's account can be understood as a further explanation of what the instantiation of mathematical concepts amounts to; then one may claim that the Kantian interpretation of Euclid's  \emph{exposition} and \emph{specification} is compatible with its interpretation as the universal instantiation in the modern sense. However the Kantian interpretation doesn't suggest by itself to interpret the instantiation as a logical procedure, i.e., as an inference of a proposition from another proposition. As the above quote makes it clear Kant describes the instantiation as a cognitive procedure of a different sort.}
 
 \paragraph{Now coming back to Euclid we must first of all admit that the  \emph{exposition} and the \emph{specification} of Theorem 1.5 as they stand are too concise for preferring one philosophical interpretation rather than another. Euclid introduces an isosceles triangle through Definition 1.20 providing no rule for constructing such a thing. (This example may serve as an evidence against the often-repeated claim that every geometrical object considered by Euclid is supposed to be constructed on the basis of Postulates beforehand.) Nevertheless given the important role of constructions in Euclid's geometry, which I explain in the next section, the idea that every geometrical object in Euclid has an associated construction rule, appears very plausible. There is also another interesting textual feature of Euclid's \emph{specification} that in my view makes the Kantian interpretation more plausible.}
 
 \paragraph{Notice the use of the first person in the  \emph{specification} of Theorem 1.5 : ``I say that ....''. In \emph{Elements} Euclid uses this expression systematically in the \emph{specification} of every theorem. Interpreting the \emph{specification} in terms of universal instantiation one should, of course, disregard this feature as merely rhetorical. However it may be taken into account through the following consideration. While the  \emph{enunciation} of a theorem is a general proposition that can be best understood \'a la Frege in the abstraction from any human or inhuman thinker, i.e., independently of any thinking \emph{subject}, who might believe this proposition, assert it, refute it, or do anything else about it, the core of Euclid's theorems (beginning with  their  \emph{exposition})  involves an individual thinker (individual subject) that cannot and should not be wholly abstracted away in this context. When Euclid  \emph{enunciates} a theorem this \emph{enunciation} does not involve - or at least is not supposed to involve - any particularities of Euclid's individual thinking; the less this \emph{enunciation} is affected by Euclid's (or anyone else's) individual writing and speaking style the better. However the \emph{exposition} and the\emph{specification} of the given theorem essentially involve an \emph{arbitrary} choice of notation (``Let $ABC$ be an isosceles triangle...''), which is an individual choice made by an individual mathematician (namely, made by Euclid on the occasion of writing his \emph{Elements}). This individual choice of notation goes on par with what we have earlier described as  \emph{instantiation}, i.e. the choice of one individual triangle (triangle $ABC$) of the given type, which serves Euclid for proving the general theorem about \emph{all} triangles of this type. The \emph{exposition} can be also naturally accompanied by drawing a diagram, which in its turn involves the choice of a particular shape (provided this shape is of the appropriate type), to leave alone the choices of its further features like color, etc.       
 }    

 \paragraph{Thus when in the \emph{specification} of Theorem 1.5 we read ``I say that the angle $ABC$ is equal to $ACB$'' we indeed do have good reason to take Euclid's wording seriously. For the sentence ``angle $ABC$ is equal to $ACB$'' unlike the sentence ``for isosceles triangles, the angles at the base are equal to one another'' has a feature that is relevant only to one particular presentation (and to one particular diagram if any), namely the use of letters $A, B, C$ rather than some others.  The words ``I say that ...'' in the given context stress this situational character of the following sentence ``angle $ABC$ is equal to $ACB$''. What matters in these words is, of course, not Euclid's personality but the reference to a particular act of speech and cognition of an individual mathematician. Proving the same theorem on a different occasion Euclid or anybody else could use other letters and another diagram of the appropriate type.}

 \footnote{
Although the choice of letters in Euclid's notation is arbitrary the  \emph{system} of this notation is not. This traditional geometrical notation has a relatively stable and rather sophisticated syntax, which I briefly describe in what follows.}   
 
  \paragraph{A competent reader of Euclid is supposed to know that the choice of letters in Euclid's notation is arbitrary and that Euclid's reasoning does not  depend of this choice. The arbitrary character of this notation should be distinguished from the general arbitrariness of linguistic symbols in natural languages. What is specific for the case of \emph{exposition} and \emph{specification} is the fact that here the arbitrary elements of reasoning (like notation) are sharply distinguished from its invariant elements. To use Kant's term we can say that behind the notion according to which the choice of Euclid's notation is arbitrary (at least at the degree that letters used in this notation are permutable) and according to which the same reasoning may work equally well with different diagrams (provided all of them belong to the same appropriate type) there is a certain invariant  \emph{schema} that sharply limits such possible choices. This schema not simply \emph{allows} for making some arbitrary choices but \emph{requires} every possible choice in the given reasoning to be wholly arbitrary. This requirement is tantamount to saying that subjective reasons behind choices made by an individual mathematician for presenting a given mathematical argument are strictly irrelevant to the ``argument itself'' (in spite of the fact that the argument cannot be formulated without making such choices). In general talks in natural languages there is no similar sharp distinction between arbitrary and invariant elements . When I write this paper I can certainly change some wordings without changing the sense of my argument but I am not in a position to describe precisely the scope of such possible changes  and identify the intended ``sense'' of my argument with a mathematical rigor. This is because my present study is philosophical and historical but not purely mathematical.}     
 
\paragraph{Thus the instantiation of a universal proposition (\emph{enunciation}) by a particular geometrical example (like triangle $ABC$)  Euclid's  \emph{exposition} serves for the formulation of this universal proposition in terms, which are suitable for a particular act of mathematical cognition made by an individual mathematician. This aspect of the \emph{exposition} is not accounted for by  the modern notion of universal instantiation itself. It may be argued that this further aspect of the \emph{exposition} needs not be addressed in a \emph{logical} analysis of Euclid's mathematics that aims at explication of the  \emph{objective meaning} of Euclid's reasoning and should not care about cognitive aspects of this reasoning. I agree that this latter issue lies out of the scope of logical analysis in the usual sense of the term but I disagree that the objective meaning of Euclid's reasoning can be made explicit without addressing this issue. Euclid's mathematical reasoning is \emph{objective} due to a mechanism that allows one to make universally valid inferences through one's individual thinking. Whatever the ``objective meaning'' might consist of this mechanism must be taken into account.}       

\section{Logical Deduction and Geometrical Production} 
\paragraph{
Remind that the \emph{proof} of Euclid's Theorem 1.5 uses not only premiss \textbf{Hyp} assumed ``by hypothesis'' but also premisses \textbf{Con 1-3} (as well as a number of other premisses of the same type) assumed ``by construction''. I turn now to the question about the role of Euclid's \emph{constructions} (which, but the way, are ubiquitous not only in geometrical but also in arithmetical Books of the \emph{Elements}) and more specifically consider the question how these \emph{constructions} support certain premisses that are used in following \emph{proofs}. 
}

As it is well-known Euclid's geometrical constructions are supposed to be realized ``by ruler and compass''. In the \emph{Elements} this condition is expressed in the \emph{Elements} through the following three 

\begin{quote}
Postulates: \\
1. Let it have been postulated to draw a straight-line from any point to any point.\\
2. And to produce a finite straight-line continuously in a straight-line.\\
3. And to draw a circle with any center and radius.
\end{quote}

(I leave out of my present discussion two further Euclid's Postulates including the controversial Fifth Postulate.) \\
\paragraph{
Before I consider popular interpretations of these Postulates and suggest my own interpretation let me briefly discuss the very term ``postulate'', which is traditionally used in English translations of Euclid's  \emph{Elements}. Fitzpatrick translates Euclid's verb ``aitein'' by English verb ``to postulate'' following the long tradition of Latin translations, where this Greek verb is translated by Latin verb ``postulare''. However according to today's standard dictionaries the modern English verb ``to postulate'' does not translate the Greek verb ``aitein'' and the the Latin verb ``postulare'' in general contexts: the modern dictionaries translate these verbs into ``to demand'' or ``to ask for''. This clearly shows that the meaning of the English verb ``to postulate'' that derives from Latin ``postulare'' changed during its lifetime.   
}

\footnote{I reproduce here Fitzpatrick's footnote about Euclid's expression ``let it be postulated'':
\begin{quote}
The Greek present perfect tense indicates a past action with present significance. Hence, the 3rd-person present perfect imperative \emph{Hitesthw} could be translated as ``let it be postulated'', in the sense ``let it stand as postulated'', but not ``let the postulate be now brought forward''. The literal translation ``let it have been postulated'' sounds awkward in English, but more accurately captures the meaning of the Greek.
\end{quote}
}  

\paragraph{
 Aristotle describes a postulate (aitema) as what ``is assumed when the learner either has no opinion on the subject or is of a contrary opinion'' (\emph{An. Post.} 76b); further he draws a contrast between postulates and  \emph{hypotheses} saying that the latter appear more plausible to the learner than the former (\emph{ibid.}). It is unnecessary for my present purpose to go any further into this semantical analysis trying to reconstruct an epistemic attitude that Euclid might have in mind ``demanding'' the reader to take his Postulates for granted. The purpose of the above philological remark is rather to warn the reader that the modern meaning of the English word ``postulate'' can easily mislead when one tries to interpret Euclid's Postulates adequately. So I suggest to read Euclid's Postulates as they stand and try to reconstruct their meaning from their context forgetting for a while what one has learned about the meaning of the term ``postulate'' from modern sources.  
}

\paragraph{Euclid's Postulates are usually interpreted as propositions of a certain type and on this basis are qualified as axioms in the modern sense of the term. There are at least two different ways of rendering Postulates in a propositional form. I shall demonstrate them at the example of Postulate 1. This Postulate can be interpreted either as the following \emph{modal} proposition:\\
(PM1): given two different points it is always possible to drawing a (segment of) straight-line between these points\\
or as the following \emph{existential} proposition: \\
(PE1): for any two different points there exists a (segment of) straight-line lying between these points. \\   
}

Propositional interpretations of Euclid's Postulates allow one to present Euclid's geometry as an axiomatic theory in the modern sense of the word and, more specifically, to present Euclid's geometrical constructions as parts of proofs of his theorems. Even before the modern notion of axiomatic theory was strictly defined in formal terms many translators and commentators of Euclid's \emph{Elements} tended to think about his theory in this way and interpreted Euclid's Postulates in the modal sense. Later a number of authors (\cite{Hintikka&Remes:1974}, \cite{Avigard&Dean&Mumma:2009}) proposed different formal reconstructions of Euclid's geometry based on the existential reading of Postulates. 
According to Hintikka \& Remes 

\begin{quote}
[R]eliance on auxiliary construction does not take us outside the axiomatic framework of geometry. Auxiliary constructions are in fact little more than ancient counterparts to applications of modern instantiation rules. 
\cite{Hintikka&Remes:1976}, p. 270
\end{quote}

\paragraph{The instantiation rules work in this context as follows. First, through the \emph{universal instantiation} (which under this reading correspond to Euclid's \emph{exposition} and \emph{specification}) one gets some propositions like \textbf{Hyp} about certain particular objects (individuals) like  $AB$ and $AC$. Then one uses Postulates 1-3 reading them as existential axioms according to which the existence of certain geometrical objects implies the existence of certain further geometrical objects, and so proves the (hypothetical) existence of such further objects of interest. Finally one uses another instantiation rule called the rule of \emph{existential instantiation}:\\ 
$$\exists x P(x) \Longrightarrow P(a) $$\\ 
and thus ``gets'' these new objects. Under this interpretation Euclid's \emph{constructions} turn into logical inferences of sort. 
As Hintikka \& Remes emphasize in their paper the principal distinctive feature of Euclid's \emph{constructions} (under their interpretation) is that these constructions introduce some \emph{new} individuals; they call such individuals ``new'' in the sense that these individuals are not (and cannot be) introduced through the universal instantiation of hypotheses making part the \emph{enunciation} of the given theorem.  
}

The propositional interpretations of Euclid's Postulates are illuminating because they allow for analyzing traditional geometrical constructions in modern logical terms. However they require a paraphrasing of Euclid's wording, which from a logical point of view is far from being innocent. In order to see this let us leave aside the epistemic attitude expressed by the verb ``postulate'' and focus on the question of \emph{what} Euclid postulates in his Postulates 1-3. Literally, he postulates the following: 

\begin{quote}
P1: to draw a straight-line from any point to any point.\\
P2: to produce a finite straight-line continuously in a straight-line.\\
P3: to draw a circle with any center and radius.
\end{quote}

As they stand expressions P1-3 don't qualify as propositions; they rather describe certain \emph{operations}! And making up a proposition from something which is not a proposition is not a innocent step. My following analysis is based on the idea that Postulates are \emph{not} primitive truths from which one may derive some further truths but are primitive operations that can be combined with each other and so produce into some further operations. In order to make my reading clear I paraphrase P1-3 as follows: 
\begin{quote}
(OP1): drawing a (segment of) straight-line between its given endpoints\\
(OP2): continuing a segment of given straight-line indefinitely (``in a straight-line)''\\
(OP3): drawing a circle by given radius (a segment of straight-line) and center (which is supposed to be one of the two endpoints of the given radius).\\
\end{quote}

Noticeably none of OP1-3 allows for producing geometrical constructions out of nothing; each of these fundamental operation produces a geometrical object out of some other objects, which are supposed to be \emph{given} in advance. The table below specifies inputs (operands) and  outputs (results) of  OP1-3:

\begin {center}
\begin{tabular}{|l|c|r|}
  \hline
  operation & input & output \\
  \hline
  OP1 & two (different) points & straight segment \\
   \hline
  OP2 & straight segment  & (bigger) straight segment \\
  \hline
  OP3 & straight segment and one of its endpoints & circle  \\ 
  \hline
\end{tabular}
\end {center}

\paragraph{PE1 as it stands does not imply that there exists at least one point or at least one line in Euclid's geometrical universe. If there are no points then there are no lines either. Similar remarks can be made about the existential interpretation of other Euclid's Postilates. Thus the existential interpretation of Postulates by itself does not turn these Postulates into existential axioms that guarantee that Euclid's universe is non-empty and contains all geometrical objects constructible by ruler and compass. To meet this purpose one also needs to postulate the existence of at least two different points - and then argue that the absence of any counterpart of such an axiom in Euclid is due to Euclid's logical incompetence. My proposed  interpretation of Postulates 1-3 as operations doesn't require such ad hoc stipulations and thus is more faithful not only to Euclid's text but also to a deeper structure of his reasoning.} 

\footnote{
Remind that the concepts of infinite straight line and infinite half-line (ray) are absent from Euclid's geometry; thus the result of OP2 is always a finite straight segment. However the result of OP2 is obviously not fully determined by its single operand. This shows that OP2 doesn't really fit the today's usual notion of algebraic operation. }

\paragraph{
Hintikka \& Remes describe Euclid's geometrical constructions as \emph{auxiliary}. Such a description may be adequate to the role of geometrical constructions in today's practice of teaching the elementary geometry but not to the role of constructions in Euclid's  \emph{Elements}. Remind that Euclid's so-called Propositions are of two types: some of them are Theorems while some other are Problems (see again the above quotation from Proclus' \emph{Commentary} ). In the \emph{Elements} Problems are at least as important as Theorems and arguably even more important: in fact  the \emph{Elements} begin and end with a Problem but not with a Theorem. As we shall now see when a given \emph{construction} makes part of a problem rather than a theorem it cannot be described as auxiliary in any appropriate sense. We shall also see the modern title ``proposition'' is not really appropriate when we talk about Euclid's Problems: while \emph{enunciations} of Theorems do qualify as propositions in the modern logical sense of the term \emph{enunciations} of Problems do not. 
}

\paragraph{I shall demonstrate these features at the well known example of Problem 1.1 that opens Euclid's \emph{Elements}; my notational conventions remain the same as in the example of Theorem 1.5.   
}

\paragraph{[\emph{enunciation:}]}
\begin{quote}
To construct an equilateral triangle on a given finite straight-line.
\end{quote}

\paragraph{[\emph{exposition}:]}

\begin{quote}
Let AB be the given finite straight-line.
\end{quote}

\includegraphics [scale=0.5]{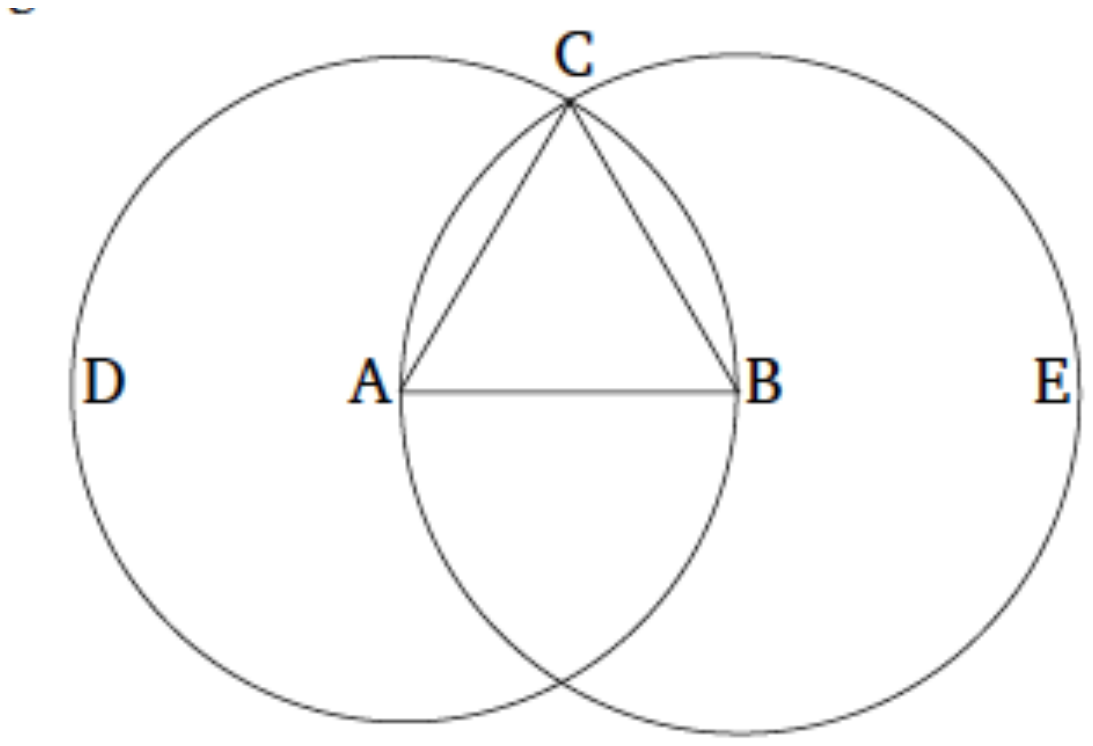}

\paragraph{[\emph{specification}:]}

\begin{quote}
So it is required to construct an equilateral triangle on the straight-line $AB$.
\end{quote}

\paragraph{[\emph{construction}:]}

\begin{quote}
Let the circle $BCD$ with center $A$ and radius $AB$ have been drawn [Post. 3], and again let the circle $ACE$ with center $B$ and radius $BA$ have been drawn [Post. 3]. And let the straight-lines $CA$ and $CB$ have been joined from the point $C$, where the circles cut one another, to the points $A$ and $B$ [Post. 1].
\end{quote}

\paragraph{[\emph{proof}:]}

\begin{quote}
And since the point $A$ is the center of the circle $CDB$, $AC$ is equal to $AB$ [Def. 1.15]. Again, since the point $B$ is the center of the circle $CAE$, $BC$ is equal to $BA$ [Def. 1.15]. But $CA$ was also shown (to be) equal to $AB$. Thus, $CA$ and $CB$ are each equal to $AB$. But things equal to the same thing are also equal to one another [Axiom 1]. Thus, CA is also equal to CB. Thus, the three (straight-lines) $CA$, $AB$, and $BC$ are equal to one another.
\end{quote}

\paragraph{[\emph{conclusion}:]}

\begin{quote}
Thus, the triangle $ABC$ is equilateral, and has been
constructed on the given finite straight-line $AB$. (Which
is) the very thing it was required to do.
\end{quote}

\paragraph{As one can see at this example \emph{enunciations} of Problems are expressed in the same grammatical form as Postulates 1-3, namely in the form of infinitive verbal expressions. I read these expressions in the same straightforward way, in which I read Postulates: as descriptions of certain geometrical \emph{operations}. The obvious difference between (\emph{enunciations} of) Problems and Postulates is this: while Postulates introduce basic operations taken for granted (drawing by ruler and compass) Problems describe complex operations, which in the last analysis reduce to these basic operations. Such reduction is made through a \emph{construction} of a given Problem: it performs the complex operation described in the \emph{enunciation} of the problem through combining basic operations OP1-3 (and possibly some earlier performed complex operations). The procedure that allows for  performing complex operations by combining a small number of repeatable basic operations I shall call a \emph{geometrical production}. In Problem 1.1 the construction of regular triangle is (geometrically) \emph{produced} from drawing the straight-line between two given points (Postulate 1) and drawing a circle by given center and radius (Postulate 3). This is, of course, just another way of saying that the regular triangle is constructed by ruler and compass; the unusual terminology helps me to describe Euclid's geometrical constructions more precisely.}

\paragraph{
Let us see in more detail how works Euclid's geometrical production. Basic operations OP1-3 like other (complex) operations need to be \emph{performed}: in order to produce an output they have to be fed by some input. This input is provided through the \emph{exposition} of the given Problem (the straight line $AB$ in the above example). OP1-3 are composed in the usual way well-known from today's algebra: outputs of earlier performed operations are used as inputs for further operations. 
}

\footnote{
Problem 1.1 involves a difficulty that has been widely discussed in the literature: Euclid does not provide any principle that may allow him to construct a point of intersection of the two circles involved into the \emph{construction} of this Problem. This flaw is usually described as a \emph{logical} flow. In my view it is more appropriate to describe this flow as properly\emph{geometrical} and fill the gap in the reasoning by the following additional postulate (rather than an additional axiom):

\begin{quote}
Let it have been postulated to produce a point of intersection of two circles with a common radius.
\end{quote}  

Even if this additional postulates is introduced here purely ad hoc, the way in which it is introduced gives an idea of how Euclid's Postulates could emerge in the real history.    
}

Just like Postulates 1-3 \emph{enunciations} of Problems can be read as modal or existential propositions (in the modern logical sense of the term). Then the (modified) \emph{enunciation} of Problem 1.1 reads:
 
 \begin{quote}
 (1.1.M) it is possible to construct a regular triangle on a given finite straight-line:
 \end{quote}
 or
\begin{quote}
(1.1.E) for all finite straight-line there exists a regular triangle on this line.
 \end{quote}

\paragraph{
As soon as the \emph{enunciations} of Euclid's Problems are rendered into the propositional form the Problems turn into theorems of a special sort. In the case of existential interpretation Problems turn into \emph{existential} theorems that state (under certain hypotheses) that there exist certain objects having certain desired properties. However this is not what we find in Euclid's text as it stands. Every Euclid's Problem ends with the formula ``the very thing it was required to do'', not ``to show'' or ``to prove''. I can see no evidence in the  \emph{Elements} that justifies the idea that in Euclid's mathematics \emph{doing} is less significant than \emph{showing} and that the former is in some sense reducible to the latter. In the Second Part of this paper I shall argue that \emph{doing} remains as much important in today's mathematics as it was in Greek mathematics, and that the idea of reducing mathematics to  \emph{showing}  or  \emph{proving} (in the precise sense of modern logic) is a unfortunate philosophical misconception.
}

According to another popular reading Euclid's Problems are tasks or questions of sort. This version of modal propositional interpretation of Euclid's Problems involves a deontic modality rather than a possibility modality: 

\begin{quote}
 (1.1.D) it is required to construct a regular triangle on a given finite straight-line:
 \end{quote}
 
\paragraph{
Indeed geometrical problems similar to Euclid's Problems can be found in today's Elementary Geometry textbooks as exercises. However the analogy between Euclid's Problems and school problems on construction by ruler and compass is not quite precise. \emph{Enunciations} of Euclid's Problems just like the \emph{enunciations} of Euclid's Theorems prima facie express no modality whatsoever. A deontic expression appears only in the \emph{exposition} of the given Problem (``it is required to construct an equilateral triangle on the straight-line $AB$''). I don't think that this fact justifies the deontic reading of the \emph{enunciation} because, as I have already explained above, the \emph{exposition} describes reasoning of an individual mathematician rather than presents this reasoning in an objective form. That every complex construction must be performed through Postulates  and earlier performed constructions is an epistemic requirement, which is on par with the requirement according to which every theorem must be proved rather than simply stated. Remind that the \emph{expositions} of Euclid's Theorems have the form ``I say that...''. This indeed makes an apparent contrast with the \emph{expositions} of Problems that have the form ``it is required to ....''. However this contrast doesn't seem me to be really sharp. Euclid's expression ``I say that...'' in the given context is interchangeable with the expression ``it is required to show that...'', which matches the closing formula of Theorems ``(this is) the very thing it was required to show''. Euclid's expression ``it is required to...'' that he uses in the \emph{expositions} of Problems similarly matches the closing formula of Problems ``(this is) the very thing it was required to do''. The requirement according to which every Theorem must be ``shown'' or ``monstrated'' doesn't imply, of course, that the \emph{enunciation} (statement) of this Theorem has a deontic meaning. The requirement according to which every Problem must be ``done'' doesn't imply either that the \emph{enunciation} of this Problem  has something to do with deontic modalities. 
}

\paragraph{
The analogy between Axioms and Theorems, on the one hand, and Postulates and Problems, on the other hand, may suggest that Euclid's geometry splits into two independent parts one of which is ruled by (proto)logical deduction while the other is ruled by geometrical production. However this doesn't happen and in fact Problems and Theorems turn to be mutually dependent elements of the same theory. The above example of Problem 1.1 and Theorem 1.5 show how the intertwining of Problems and Theorems works. Theorems, generally, involve \emph{constructions} (called in this case auxiliary), which may depend (in the order of geometrical production) on earlier treated Problems (as the \emph{construction} of Theorem 1.5 depends on Problem 1.3.) Problems in their turn always involve appropriate \emph{proofs} that prove that the \emph{construction} of the given theorem indeed performs the operation described in the \emph{enunciation} of this theorem (rather than performs some other operation). Such \emph{proofs}, generally, depend (in the order of the protological deduction) on certain earlier treated Theorems (just like in the case of \emph{proofs} of Theorems).  
Although this mechanism linking Problems with Theorems may look unproblematic it gives rise to the following puzzle. Geometrical production produces geometrical objects from some other objects. Protological deduction deduces certain propositions from some other propositions. How it then may happen that the geometrical production has an impact on the protological deduction? In particular, how the geometrical production may justify premises assumed ``by construction'', so these premises are used in following \emph{proofs}? 
}

\paragraph{In order to answer this question let's come back to the premise \textbf{Con3} ($AF = AG$) from Theorem 1.5 and see what if anything makes it true. $AF = AG$ because Euclid or anybody else following Euclid's instructions constructs this pair of straight segments in this way. How do we know that by following these instructions one indeed gets the desired result? This is because the \emph{construction} of Problem 1.3 that contains the appropriate instruction is followed by a \emph{proof} that proves that this \emph{construction} does exactly what it is required to do. \emph{Construction} 1.3 in its turn uses \emph{construction} 1.2  while \emph{construction} 1.2 uses \emph{construction} 1.1 quoted above. In other words \emph{construction}  1.1 (geometrically) produces \emph{construction}  1.2 and \emph{construction}  1.2 in its turn produces \emph{construction} 1.3. This  geometrical production produces the relevant part of \emph{construction} 1.5 (the construction of equal straight segments $AF$ and $AG$) from first principles, namely from Postulates 1-3. In order to get the corresponding protological deduction of  premise \textbf{Con3}  from first principles we should now look at \emph{proofs} 1.1, 1.2 and 1.3 and then combine these three proofs into a single chain. For economizing space I leave now details to reader and only report what we get in the end. The result is somewhat surprising from the point of view of the modern logical analysis. The chain of \emph{constructions} leading to \emph{construction} 1.5 involves a number of circles (through Postulate 3). Radii of a given circle are equal by definition (Definition 1.15). Thus by constructing a circle and its two radii, say, $X$ and $Y$ one gets a primitive (not supposed to be proved) premise $X = Y$. Having at hand a number of premises of this form and using Axioms as inference rules (but not as premises!) one gets the desired deduction of \textbf{Con3}. The fact that first principles of the protological deduction of \textbf{Con3} appear to be partly provided by a Definition helps to explain why Euclid places his Definitions among other first principles such as Postulates and Axioms.      
}

\paragraph{The above analysis allows for disentangling the protological deduction of \textbf{Con3} from the geometrical production of straight segments  $AF$, $AG$ and so the aforementioned puzzle remains even after we have looked at Euclid's reasoning under a microscope. Even if we can describe in detail the impact of Problems to Theorems and vice versa it remains unclear how the two kind of things can possibly work together. Here is my tentative answer to this question. Every Euclid's \emph{proof}   involves only \underline{concrete} premises like \textbf{Con3} and \textbf{Hyp}, which are related to certain individual objects. It is assumed that such a premise is valid \underline{only if} the corresponding object is effectively constructed. (At least this concerns all premises ``by construction''; as we have seen at the example of Theorem 1.5 hypothetic premises sometimes don't respect this rule.) This fundamental principle links Euclid's geometrical production and protological deduction together.}

\paragraph{
One may argue that my proposed analysis after all is not significantly different from the standard logical analysis of Euclid's geometrical reasoning according to which Euclid first proves that certain geometrical objects exist and only then prove some further propositions concerning properties of these objects. Is there indeed any significant difference between proving that  such-and-such object exist and producing this object through what I call the geometrical production? There is of course a difference of a metaphysical sort between these two notions: to produce an object is not quite the same thing as to prove that certain object exists. But arguably this difference has no objective significance and so one may simply ignore it. There is however a further difference between the geometrical production and the mathematical existence, which seems me more important. Euclid's \emph{Elements} contain two sets of rules, namely Axioms and Postulates, supposed to be applied to operations of two different sorts: Axioms tell us how to derive equalities from other equalities while Postulates tell us how to produce geometrical objects from other geometrical objects. A logical analysis of Euclid's geometry that involves a propositional (in particular existential) reading of Postulates aims at replacing these two sets of rules by a single set of rules called \emph{logical}. I would like to stress again that the results my proposed analysis don't exclude the possibility of logical analysis. Such a replacement may be or be not a good idea but in any event logical rules are not made in the Euclid's text explicit and I don't see much point in saying that he uses rules of this sort implicitly. The fact that \emph{we} can use today modern logic for interpreting Euclid is a completely different issue.  An interpretation of Euclid's geometry by means of logical analysis can be indeed illuminating but one should not confuse oneself by believing that Euclid already had a grasp of modern logic even if could not formulate principles of this logic explicitly. At the same the combination of protological deduction and geometrical production found in Euclid may shed some light on modern logic and its involved relationships with mathematics as we shall shortly see. 
}

\paragraph{For further references I shall call the 6-part structure of Euclid's Problems and Theorems \emph{Euclidean structure}. As the above analysis makes clear the Euclidean structure doesn't work with the modern notion of axiomatic theory but requires a different setting, which combines the protological deduction with the geometrical production. Let us now see what the Euclidean structure can tell us about today's mathematics. 
}

\part{Modern Axiomatic Method}

\section{Euclid and Modern Mathematics}
\paragraph{What has been said above may give one an impression that in Euclid's  \emph{Elements} we deal with an archaic pattern of mathematical thinking that has noting to do with today's mathematics. However this impression is wrong. In fact the Euclidean structure is apparently present in today's mathematics, perhaps in a slightly modified form. Consider the following example taken from a standard mathematical textbook (\cite{Kolmogorov&Fomin:1976}, p. 100, my translation into English):
}
\paragraph{Theorem 3:}

\begin{quote}
Any closed subset of a compact space is compact
\end{quote}

\paragraph{Proof:}

\begin{quote}
Let $F$ be a closed subset of compact space $T$ and $\{F_{\alpha}\}$ be an arbitrary centered system of closed subsets of subspace $F \subset T$. Then every $F_{\alpha}$ is also closed in $T$, and hence $\{F_{\alpha}\}$ is a centered system of closed sets in $T$. Therefore $\cap F_{\alpha} \neq \emptyset$. By Theorem 1 it follows that $F$ is compact.
\end{quote}

\paragraph{
Although the above theorem is presented in the usual for today's mathematics form ``proposition-proof'', its Euclidean structure can be made explicit without re-interpretations and paraphrasing:}

\paragraph{[\emph{enunciation:}]}
\begin{quote}
Any closed subset of a compact space is compact
\end{quote}

\paragraph{[\emph{exposition:}]}

\begin{quote}
Let $F$ be a closed subset of compact space $T$
\end{quote}

\paragraph{[\emph{specification}: absent]}

\paragraph{[\emph{construction}:]}

\begin{quote}
[Let] $\{F_{\alpha}\}$ [be] an arbitrary centered system of closed subsets of subspace $F \subset T$. 
\end{quote}

\paragraph{[\emph{proof}:]}

\begin{quote}
[E]very $F_{\alpha}$ is also closed in $T$, and hence $\{F_{\alpha}\}$ is a centered system of closed sets in $T$. Therefore $\cap F_{\alpha} \neq \emptyset$. By Theorem 1 it follows that $F$ is compact.
\end{quote}

\paragraph{[\emph{conclusion}: absent ]}

\paragraph{The absent \emph{specification} can be formulated as follows:}

\begin{quote}
I say that $F$ is a compact space
\end{quote}

\paragraph{while the absent \emph{conclusion} is supposed to be a literal repetition of the \emph{enunciation} of this theorem. Clearly these latter elements can be dropped for parsimony reason. In order to better separate the \emph{construction} and the \emph{proof} of the above theorem the authors could first construct set 
 $\cap F_{\alpha}$ and only then prove that it is non-empty. However this variation of the classical Euclidean scheme also seems me negligible. I propose the reader to check it at other modern examples that the Euclidean structure remains today at work.} 
 
 \paragraph{
 Does this mean that the modern notion of axiomatic theory is inadequate to today's mathematical practice just like it is inadequate to Euclid's mathematics? Such a conclusion would be too hasty. Arguably, in spite of the fact that today's mathematics preserves some traditional outlook it is essentially different. So the ``Euclidean appearance'' of today's mathematics cannot be a sufficient evidence for the claim the the Euclidean structure remains significant in it. In order to justify this claim a different argument is needed. 
} 
 
\paragraph{
Before I try to provide such an argument I would like to point to the fact that the modern notion of axiomatic theory is used in today's mathematics in two rather different ways. First, it is used as a broad methodological idea that determines the general architecture of a theory but has no impact on details. Such an application of the modern axiomatic method is usually called \emph{informal}. Second, the notion of axiomatic theory is used for building \emph{formal} theories that contain a list of axioms and a set of theorems derived from these axioms according to explicitly specified rules of logical inference. In the next Section I shall describe the notion of formal axiomatic theory more precisely and try to explain why it is called ``formal''. Here I would like only to stress that only in a formal setting the modern notion of axiomatic theory is made explicit, and that the informal application of this notion in the current mathematical practice depends on the assumption that the informal mathematics can be appropriately formalized. Thus in order to test the modern notion of axiomatic theory against today's mathematical practice we should first of all look at formal theories and analyze what the formalization of the informal mathematics amounts to.    
}

\section{Formalization}
\paragraph{While the concept of formalization is rather epistemological than purely mathematical the concept of \emph{formal axiomatic theory} admits a rigorous mathematical definition. This rigor comes with a price: mathematical definitions of formal theory vary from one textbook to another, see for example \cite{Kleene:2009} and \cite{Mendelson:1997}. Yet there is a common understanding that such special definitions share a common conceptual core, which I shall try to present here. The definition of formal axiomatic theory that I give below is not the most general and not the most precise; strictly speaking it doesn't qualify as a mathematical definition. I consider only the case of \emph{deductive} theory that is relevant to the axiomatic method. 
}

Thus a (deductive) formal theory consists of the following:

\begin{itemize}
\item (a) a formal \emph{language} that includes an \emph{alphabet} of symbols and rules of building \emph{formulae} from symbols of this alphabet;
\item (b) a list of \emph{axioms} of the given theory, which is a distinguished set of well-formed formulae of the theory; formulae which are not axioms are called \emph{theorems};
\item (c) \emph{rules of inference} of formulae from other formulae; the inference of formulae from other formulae is called \emph{deduction}. As far as rules of inference work universally for all theories or at least for a large class of theories the deduction is called \emph{logical}. 

\end{itemize}

\paragraph{
One usually further assumes that a given formal theory can be \emph{interpreted}. An interpretation assigns truth-values to formulae of the theory. Interpretations under which the axioms of the given theory are true is called a \emph{model} of this theory. As far as rules of inference are truth-preserving (which is the usual requirement) all theorems deduced from axioms are also true in every model of the given theory.  Thus a formal theory such that all its theorems are deduced from axioms is a concrete mathematical realization of the informal notion of axiomatic theory discussed earlier.  Canonical examples of formal theories are Peano Arithmetics (PA) and an axiomatic theory of sets called ZF after Zermelo and Fraenkel.
}

\footnote{There exist many different versions of PA and ZF but I don't need to be more specific about them here.} 

\paragraph{
The notion of formal theory allows for more precise definitions of concepts of theory, axiom and theorem than usual. However it is appropriate to ask whether or not such precise definitions adequately capture the meaning of the corresponding terms, which is known from the usual ``informal'' mathematical practice. Is a formal theory ``really'' a theory? Is a formal mathematical theorem in the sense of the above definition ``really'' a mathematical theorem? Asking these questions I don't  assume that the meanings of terms ``theory'' and ``theorem'' live somewhere on Platonic Heavens or in the Current Practice and need only be described correctly. I rather assume that our concepts of theory, theorem, etc. are \emph{formed} through our attempts to improve upon the existing theoretical practice in the pure mathematics and elsewhere.  So the above question can be more precisely formulated as follows: Are formal mathematical theories indeed able to replace the so-called informal theories in  practice? Do they really do this at least at certain degree? Is there an advantage of doing mathematics formally rather than informally? If the answers to these questions are in positive then one may argue that even if the existing mathematical practice remains largely informal a \emph{formalization} of mathematics remains a sound project for the future. One may also argue that in spite of the fact that today's mathematical theories still look like the theory of Euclid's \emph{Elements} in certain respects the core structure of these today's theories is very different and can be better described with the modern notion of axiomatic theory. Finally one may argue that to make this core structure explicit one has to formalize today's mathematical theories, i.e., rewrite them as formal theories in the above precise sense.   
} 
\paragraph{
Before I try to answer the above questions and evaluate the above argument on this basis I would like to make some historical remarks concerning the notion of being formal. I shall talk first about formal logic and then about formal approaches in mathematics. The term ``formal logic'' has a long history that I shall not try to trace here but only stress the fact that the meaning of this term changed drastically during last several centuries. For example, Kant's notion of formal logic (as distinguished from the transcendental logic) is quite different from our modern notion of formal logic (as distinguished from the so-called informal logic). Let me now talk only about the modern notion. Today people tend to qualify as formal any logical theory that uses mathematical symbolic methods and, correspondingly, qualify as informal any account of logic that does not systematically applies such methods but presents itself with a prose written in some natural language. Although the application of mathematical methods in logic is not a recent idea new applications of mathematical methods in logic in the late 19th - early 20th centuries produced a genuine revolution in logic and brought about the formal logic in the modern sense of the term. This historical remark makes clear the reason \emph{why} the modern formal logic is more rigorous than the traditional logic. And it also allows one to specify the claim. The modern formal logic is more rigorous than the traditional logic because it systematically applies mathematical methods. For the same reason Newtonian physics is more rigorous than Wolfian physics, which expresses itself with a prose. One should not forget that we are talking here about a \emph{mathematical} rigor, which, arguably, is not the only kind of rigor that one can possibly attain in a given field of research. 
}

\paragraph{
Although the notion of formal mathematical theory has a much shorter history the meaning of the term also changed during the 20th century. In the beginning of this century this term referred to the then-novel sort of theories a canonical example of which was given by Hilbert in his book on foundations of geometry \cite{Hilbert:1899} first published in 1899. Such theories were called formal in the contradistinction with traditional \emph{contentual} theories because they involved the associated notion of interpretation described in the beginning of this Section. Although formal theories in the modern sense of the term also allow for interpretations this feature is no longer seen as sufficient for calling a theory formal. In particular, the geometrical theory developed by Hilbert in \cite{Hilbert:1899} does not qualify as formal in today's sense.    
}

\paragraph{
Now we are better prepared to discuss formalization of mathematics and evaluate its epistemic consequences. Formalization of an existing informal mathematical theory $T$ is a way of reformulating $T$ in new terms; such a procedure brings about a formal theory $F$, which in some sense expresses the same mathematical content as $T$ but arguably has some epistemic advantages over $T$.  It is worth to note that reformulations of older theories in new terms occurred in mathematics permanently through its long history. What we call Euclidean geometry today differs strikingly from what one can read in the geometrical chapters of Euclid's \emph{Elements} - and still there is a sense in which the two things represent the same mathematical theory. To make this sense precise and specify some general principles allowing for an adequate translation of mathematical theories from one setting to another remains, in my view, a challenging and wholly open problem in philosophy of mathematics. I am not trying to solve this problem here but nevertheless mention it because I think that the long existing practice of reformulating mathematical theories in new terms is an appropriate historical background for understanding formalization.    
}  

\paragraph{
Now comparing formalization with similar shifts occurred in mathematics in the past (like translation of traditional geometrical problems into an algebraic language in the late 17th century) one can observe that so far formalization was not quite successful.  Although a study of formal theories developed into a well-established area of mathematics (see below) the mainstream mathematics remained largely informal; none of recent significant advances in mathematics (like the recent proof of Poincar\'e conjecture) involved formal methods (in the sense of ``formal'' specified above). True, an evaluation of success of formalization of mathematics depends on what one expects to achieve with it, and there is no stable consensus about this issue even among those people who think about formalization enthusiastically. Anyway the very fact that studies in formal mathematics is currently a relatively isolated area largely ignored by the mainstream deserves, in my view, a very serious consideration.        
}
 
 \paragraph{
 An apparent exception is the Set theory. After the discovery of famous paradoxes of so-called ``naive'' Set theory Zermelo famously proposed a sketch of formal axiomatic theory of sets aimed at saving the Set theory from inconsistencies.  Since then the mainstream research in Set theory focused on studies on various formal theories of sets of models of such theories. In the rest of mathematics people continued to use sets widely in the informal manner with an official proviso according to which the informal concept of set is nothing but a naive version of the set concepts treated rigorously by the set-theorists with formal means. Although the research into formal set theories brought about a lot of mathematical results \emph{about these theories} it is less clear that it augmented our knowledge about \emph{sets themselves}. Think of Continuum Hypothesis (CH). G\"odel \cite{Godel:1938} discovered that ZF is consistent with CH (by building a model of ZF in which CH holds) and later Cohen \cite{Cohen:1963}  discovered that ZF is consistent with the negation of CH (by building a model of ZF in which CH does not hold). So we know now that neither CH nor its negation can be derived from the axioms of ZF. However important this result may be it does not provide any definite answer to the original question, neither it allows to claim that the original question is ill-posed, so that no definite answer exists. An additional axiom - or some wholly new system of axioms for Set theory - may eventually help, of course, to settle the problem in the sense that CH or its negation can be deduced from the new system of axiom. There are obvious trivial ``solutions'' of this sort. Then however it remains to show that the system of axiom for Set theory solving the CH problem is a ``right'' one, and so the proposed solution is ``genuine''. I cannot see how this can be done on purely mathematical grounds; any possible argument to the effect that one system of axioms for Set theory is ``more natural'' than some other has a speculative nature and lacks any objective validity. Even if one gets some non-trivial proof of CH from some system of axioms that appear to be in some sense natural one can hardly claim that this system of axioms is the ``right one'' solely because it solves the CH problem and because such a proposed solution is smart and elegant.   
 }
\paragraph{ 
It may be argued that the formal axiomatic framework makes explicit a relativistic nature of mathematics, which we should learn to live with; according to this viewpoint it is pointless to ask whether CH is true or false without further qualifications, and all that mathematicians can do is to study which axioms do imply CH (modulo some specified rules of inference), which imply its negation, and which do neither (like the axioms of ZF). More generally, the only thing that mathematics can do is to provide true propositions of the  \emph{if - then} form:  \emph{if} such-and-such propositions are true  \emph{then} certain other propositions are also true. I cannot see how such a deductive relativism (or ``if-thenism'') about mathematics can be sustainable. It is incompatible not only with the common mathematical practice but also, more specifically, with the current practice of studying formal axiomatic systems. Consider the statement S of independence of CH from the axioms of ZF that I have mentioned earlier. S is commonly seen as an established theorem on a par with any other firmly established mathematical theorem. However S is not of the  \emph{if - then} form. Moreover the proof of S (that comprises the construction of G\"odel's  model $L$ verifying CH and Cohen's forcing construction falsifying CH) is a piece of rather sophisticated ``usual'' or ``naive'' mathematics but not a formal inference within certain axiomatic theory. A consistent if-thenist would not hold without further qualifications that CH is independent from the axioms of ZF  but rather say that it depends of one's assumptions.  
} 

\paragraph{ 
Thus even in Set theory informal proofs like the proof of independence of CH from axioms of ZF apparently remain indispensable. In spite of the fact that that the modern Set theory no longer considers sets naively but works instead with various formal axiomatic theories of sets this modern theory like any other modern mathematical theory relies on non-formalized proofs. What is specific for the modern Set theory is its  \emph{object} rather than its method. Instead of studying sets directly in the same way in which, say, group-theorists study groups, set-theorists study formal axiomatic theories of sets. However the \emph{methods} used by modern set-theorists are not essentially different from methods used in other parts of today's mathematics. It remains in my sense an open methodological question whether or not such a roundabout way of studying sets has indeed proved effective. True, at the present there is no clear alternative to it. However it is not inconceivable that in the near future the mathematical community may bring about an improved ``naive'' concept of set that would allow one to study sets like groups. It is not inconceivable that such an old-fashioned way of thinking about sets could after all allow for a real progress in the CH problem. In any event it seems me important to keep such a possibility open and not try to take it out of the table using some philosophical arguments. The description of sets as homotopy zero-types suggested by J.-P. Marquis is a tentative realization of this idea.       
}

\footnote{Presented at the workshop \emph{Sets within Geometry}, Nancy (France), July 26-29, 2011}

\paragraph{
The distinction between a \emph{theory} and \emph{meta-theory} (and, more generally, between \emph{mathematics} and \emph{meta-mathematics}, which dates back to Hilbert, is helpful for making things clearer. In modern Set theory a \emph{theory} is ZF or another formal axiomatic theory while proofs of independence of CH from the axioms of ZF and similar results belong to a \emph{meta-theory} that tell us important things about formal axiomatic theories. I would like, however, to point here to the fact that this standard terminology is heavily philosophically-laden and suggests a very particular view on mathematics. Namely, it suggests the view according to which formal axiomatic theories are ``usual'' mathematical theories while meta-theories belong to a special domain of meta-mathematics that lays beyond the usual mathematics and has some philosophical flavor. But if we leave philosophy aside and describe the same situation in the context of the current mathematical practice we observe quite the opposite. Formal axiomatic theories are \emph{not} mathematical theories in the usual sense of the word while their corresponding meta-theories look like ``usual'' theories from any other area of mathematics. By the analogy with the distinction between an \emph{object language} and a \emph{meta-language} in formal semantic it would be more appropriate to use the term ``object-theory'' for what in formal axiomatic studies (but not in the rest of mathematics!) is usually called simply a ``theory''.  
}

\paragraph{As far as we are talking about meta-theories in a mathematical context such theories are mathematical theories at the first place - whatever philosophical meaning one may wish to attach to them. And their corresponding object-theories are mathematical objects (rather than theories) at the first place. This remark makes it clear where the power of formalization comes from. Not surprisingly the case of formalization of mathematics turns to be similar to the case of formalization of logic. Formalization of mathematical theories turns these theories into mathematical objects of sort. This allows for studying \emph{mathematically} a class of important problems that earlier were studied only speculatively or ``informally''. These are so-called ``meta-mathematical'' problems concerning the mutual (in)dependence of mathematical theories, provability, etc. During 80 years of such studies there were obtained a number of firm mathematical results including, in particular, the independence results about CH mentioned above. The properly mathematical aspect of these result should be sharply distinguished from associated philosophical claims. The objective validity of the former doesn't extend to the latter, so the idea that meta-mathematics is in a position to prove philosophical claims by mathematical methods cannot be justified. \emph{Qua} mathematical results all meta-mathematical results concern only formalized theories. Which consequences if any these (meta-)mathematical results may have for mathematics in general is a further question that is not mathematical in character. Only if one assumes that formal (object-) theories represent ``usual'' mathematical theories \emph{adequately} one may claim that the meta-mathematical results have a direct bearing on these ``usual'' theories. Let me stress once again that the aforementioned assumption is not mathematical in character. It is a very general assumption \emph{about} mathematics that tell us what mathematics is and what it should be. One cannot possibly prove or disprove such an assumption by using purely mathematical arguments.  
}

\paragraph{
There remains then the following three possibilities, which are not mutually exclusive. One may come up with an idea what mathematics \emph{should} be on purely speculative grounds. One may also look attentively what mathematics actually \emph{is} and what it used to be in the past and form one's general ideas about mathematics on this ground. Finally one my try to change the current mathematical practice according to one's view on mathematics developed on speculative grounds. In reality every researcher in the field cannot avoid combining all the three approaches in one way or another. The main virtue of speculation is that it allows for changing usual ideas about mathematics. However unless such new ideas are tested against the practice they remain ``merely speculative'' and their value remain uncertain. The main virtue of the careful analysis of the past and present mathematical practice is that it provides a necessary background for any attempted reform of this practice. It is clear that formalization of logic and mathematics has opened some new fields of mathematical studies but it remains so far unclear what bearing if any the results of these studies have on the rest of mathematics. The persisting gap between the common mathematical practice and its formalized counterpart that is apparent not only in formal reconstructions of older mathematics but also in today's mathematics suggests that something goes wrong with the formalization. However impressive proven mathematical facts about formal object-theories may be these facts don't provide by themselves any substance to the claim that the theory-objects can substitute usual mathematical theories.  Outside the context of their informal meta-theories formal object-theories  have no objective mathematical content just like Euclidean triangles have no objective mathematical content outside the context of Euclidean geometry.
}

\paragraph{ One may argue as follows. True, any formal (object-)theory requires some supporting meta-theory. Such a meta-theory can be also formalized in its turn. This starts an infinite regress that leads nowhere. One cannot and shouldn't formalize everything. This is why a reasonable attitude is to build the core mathematics formally using for this end some necessary informal meta-mathematical arguments. Even if this informal instrument cannot be then wholly taken away like the Wittgenstein's ladder it can be viewed as a part of general philosophical underpinning of mathematics rather than as a part of mathematics proper.  
} 

\paragraph{ 
Once again I insist that the above reasoning is based on an a priori view about mathematics that is not confirmed neither by the old nor by the new mathematics. The very fact that a mathematical  \emph{object} of a particular sort is called a ``theory'' and is supposed to represent a theory doesn't provide any substance to the claim that it indeed represents a theory in some strong sense of the word. 
}

\section {Objectivity and Objecthood in Modern Mathematics}
\paragraph{ 
The above critical remarks about formalization don't imply that formalization is a bad idea. However they show that the usual way of formalizing mathematical theories involves a systematic difficulty. I claim that what is wrong in it is the usual (informal) notion of axiomatic theory, namely, the notion according to which a theory in general (and mathematical theory in particular) is a bunch of propositions generated by a finite list of propositions called axioms according certain rules of (logical) inference. 
}
\paragraph{ 
In fact this idea dates back to Aristotle and thus is far from being modern. This Aristotelian notion of theory had apparently little or no influence on Greek mathematics (that followed Euclid rather than Aristotle) but it became quite influential in the medieval Scholasticism, in particular, in the scholastic physics. The Early Modern mathematically-laden science that triumphed with Newton's  \emph{Principia} largely rejected the old scholastic pattern of theory-building and developed a very different notion of scientific theory that was best described in general terms by Kant in his \emph{Critique of the Pure Reason}. Kantian philosophy of science and mathematics remained the mainstream until the beginning of the 20th century when the old scholastic pattern of theory-building kicked back under the new name of modern axiomatic method. Let me briefly sketch this latter  development.  
}

In Kant's view the \emph{objectivity} of pure mathematics (which underlies the objectivity of the mathematically-laden empirical science) has its source in its \emph{objecthood}, i.e., in the universal schemata according to which one constructs mathematical \emph{objects} - but not just in the universal character of the involved concepts. This, according to Kant, is a crucial difference between the mathematical reasoning and the philosophical speculation:  

\begin{quote}
``Give a philosopher the concept of triangle and let him try to find out in his way how the sum of its angles might be related to a right angle. He has nothing but the concept of figure enclosed by three straight lines, and in it the concept of equally many angles. Now he may reflect on his concept as long as he wants, yet he will never produce anything new. He can analyze and make distinct the concept of a straight line, or of an angle, or of the number three, but he will not come upon any other properties that do not already lie in these concepts.   But now let the geometer take up this question. He begins at once to \underline{construct a triangle}.  Since he knows that two right angles together are exactly equal to all of the adjacent angles that can be drawn at one point on a straight line, he extends one side of his triangle and obtains two adjacent angles that together are equal to the two right ones. .... In such a way through a chain of inferences that is always \underline{guided by intuition}, he arrives at a fully illuminated and at the same time general solution of the question.'' (\emph{Critique of Pure Reason}, A 716 / B 744)
\end{quote}

 \paragraph{
 Kant's philosophy of mathematics and mathematically-laden science is based upon the best contemporary science represented by Newton's  \emph{Principia} \cite{Friedman:1992}. By this I don't mean, of course, that Kant derives his philosophical claims from the principles of Newtonian physics; Kant's  \emph{critical} philosophy rather aims at explaining how the type of knowledge best represented by the Newtonian physics is possible (as an objectively valid knowledge). Anyway this makes Kant's philosophy strongly dependent of the contemporary mathematics and science. H. Cohen, P. Natorp and other neo-Kantians who wished to sustain the Kantian project of critical philosophy in the 19th century realized this fact very clearly and made efforts to incorporate into the Kantian philosophy a historical dimension and make it to keep track of the current progress of science (including the pure mathematics) \cite{Heis:2007}. It was not quite clear in the 19th century and it still remains a controversial question today which (if any) features of Kant's original approach remain sustainable in the context of the current science and mathematics, and which features of this original approach are hopelessly outdated. More radically one may wonder if there is anything at all in Kant's analysis that has survived all the dramatic changes in science and pure mathematics that have happened since Kant's own time.   
}
 
In spite of a number of interesting attempts of upgrading the Kantian philosophy of mathematics in order to account for new mathematical developments such as the invention of non-Euclidean geometries at certain point the Kantian line in the philosophy of mathematics has been largely abandoned. Bertrand Russell's intellectual develpment is representative in this sense: after publishing in 1897 his \emph{Essay on Foundations of Geometry} \cite{Russell:1897}, which is an explicit attempt to modernize Kant's views on this subject, already in 1903 Russell publishes his \emph{Principles of Mathematics} \cite{Russell:1903} where the author develops the subject on wholly new grounds. In the \emph{Introduction} to this book Russell explains his attitude to the Kantian line of thought as follows:

\begin{quote}
It seemed plain that mathematics consists of deductions,
and yet the orthodox accounts of deduction were largely or wholly
inapplicable to existing mathematics. Not only the Aristotelian
syllogistic theory, but also the modem doctrines of Symbolic Logic,
were either theoretically inadequate to mathematical reasoning, or at
any rate required such artificial forms of statement that they could not
be practically applied. In this fact lay the strength of the Kantian
view, which asserted that mathematical reasoning is not strictly formal,
but always uses intuitions, i.e. the a priori knowledge of space and
time. Thanks to the progress of Symbolic Logic, especially as treated
by Professor Peano, this part of the Kantian philosophy is now capable
of a final and irrevocable refutation. By the help of ten principles of deduction and ten other premisses of a general logical nature (e.g. implication is a relation"), all mathematics can be strictly and formally deduced. ...\\
The general doctrine that all mathematics is deduction by
logical principles from logical principles was strongly advocated by
Leibniz... But owing partly to a faulty logic, partly to belief in the logical necessity
of Euclidean Geometry, he was led into hopeless errors in the endeavour
to carry out in detail a view which, in its general outline, is now known
to be correct. The actual propositions of Euclid, for example, do not
follow from the principles of logic alone ; and the perception of this fact
led Kant to his innovations in the theory of knowledge. But since
the growth of non-Euclidean Geometry, it has appeared that pure
mathematics has no concern with the question whether the axioms
and propositions of Euclid hold of actual space or not .....  What pure mathematics asserts is merely that the Euclidean propositions follow from the Euclidean axioms, i.e.,
it asserts an implication. ....  We assert always in mathematics
that if a certain assertion $p$ is true of any entity $x$ or of any set of
entities $x, y, z ...$, then some other assertion $q$ is true of those entities ;
but we do not assert either $p$ or $q$ separately of our entities.
\end{quote}

 \paragraph{
The above argument, which is supposed to refute Kant, obviously begs the question. From the outset Russell takes it for granted that ``mathematics consists of deductions'' and his following remarks make it clear that by deduction Russell means here a \emph{formal logical deduction}, i.e. a deduction of propositions from certain other propositions according to some general rules, which are not specific for mathematics. This assumption overtly contradicts what Kant says about mathematics, so the following argument only adds to this statement some additional details but doesn't justify it. Kant's crucial objection to the Leibnizian view on mathematics, to which Russell adheres here, is this. From a \emph{formal} point of view (i.e. as far as only \emph{logical form} of sentences is taken into consideration) mathematics is no different from a mere speculation; a speculative theory can be developed on an axiomatic basis just like a mathematical theory. What makes the difference between mathematics and speculation is the fact that mathematics constructs its \emph{objects} according to certain rules while speculation proceeds with concepts without being involved in any similar constructive activity. The fact that the speculative thought may also posit some entities falling under these concepts from the Kantian viewpoint does not constitute an objection because such stipulated entities doesn't qualify as \emph{objects} in the Kantian sense of the term. Behind an \emph{object} there is a procedure (governed by a certain \emph{rule} that constructs this object while speculative entities are stipulated as mere things falling under some given description without any constructive procedure. This is the reason why the pure mathematics is \emph{objective} in the sense in which the pure speculation is not. What makes the pure mathematics objective is the rule-like character of object-construction. The formal logical consistency is a necessary but not sufficient condition for claiming that a given axiomatic theory is objectively valid.   
}

Russell's critique of Kant  in the  \emph{Principles of Mathematics} doesn't take into the account the Kantian problem of separation of the pure mathematics from the pure speculation. In this sense his approach is more traditional than Kant's and from a Kantian viewpoint qualifies as dogmatic. Like Aristotle and Leibniz Russell provides his philosophy of mathematics with a metaphysical doctrine that he calls the  \emph{logical atomism}. This is how he describes the relation of this doctrine to logic and mathematics in the  \emph{Introduction} to his \cite{Russell:1918}:

\begin{quote}
As I have attempted to prove in  \emph{The Principles of Mathematics},
when we analyse mathematics we bring it all back to logic. It all
comes back to logic in the strictest and most formal sense. In the
present lectures, I shall try to set forth in a sort of outline, rather
briefly and rather unsatisfactorily, a kind of logical doctrine
which seems to me to result from the philosophy of mathematics - 
not exactly logically, but as what emerges as one reflects: a
certain kind of logical doctrine, and on the basis of this a certain
kind of metaphysic. 
 \end{quote}

Thus Russell sees Kant's work in the philosophy of mathematics as an attempt to fill logical gaps appearing when one tries to reconstruct Euclid's geometry with Aristotle's syllogistic logic (which hardly correctly describes Kant's own intention). Russell suggests two independent reasons why there are such gaps: first, because Euclid's geometry is logically imperfect and, second, because Aristotle's logic is not appropriate for doing mathematics. However the new mathematics (including non-Euclidean geometries) and the new symbolic logic taken together, according to Russell, wholly fix the problem making Russell's Leibnizian dream real. What Russell's \emph{The Principles of Mathematics} aim at is made clear by the following lines that I take from the \emph{Preface} to this work:

 \begin{quote}
The second volume, in which I have had the great good fortune
to secure the collaboration of Mr A. N. Whitehead, will be addressed
exclusively to mathematicians; it will contain chains of deductions,
from the premisses of symbolic logic through Arithmetic, finite and
infinite, to Geometry, in an order similar to that adopted in the present
volume ; it will also contain various original developments, in which the
method of Professor Peano, as supplemented by the Logic of Relations,
has shown itself a powerful instrument of mathematical investigation.   
 \end{quote}
 
 (The planned second volume of \emph{The Principles of Mathematics} appeared later as a co-authored independent three-volume work \cite{Russell&Whitehead:1910-1913}.) 
 
 \footnote{Cf. \cite{Hylton:1990}, p. 7-8 
 \begin{quote} 
 From August 1900 until the completion of Principia Mathematica in 1910 Russell was both a metaphysician and a working logician. The two are completely intertwined in his work: metaphysics was to provide the basis for logic; logic and logicism were to be the basis for arguments for the metaphysics.
 \end{quote}
 }
 
People who tried to push the Kantian line without loosing its adequacy to the contemporary mathematics and without turning the Kantian philosophy into an orthodoxy had indeed a hard time in the beginning of the 20th century.  Ernest Cassirer published in 1907 a paper \cite{Cassirer:1907} on the issue of relationships between Kantian philosophy and the ``new mathematics''. Referring to Russell's 1903 view and new formal logical methods under the name of ``logistics'' Cassirer says: 

 \begin{quote}
 With this arises a problem, which lies completely outside the scope of
logistics ...Worrying about the rules that govern the world of objects is
completely left to direct observation, which is the only one that can
teach us ... whether we can find here certain regularities or a pure chaos. Logic and
mathematics deal only with the order of concepts; they don't
contest the order or the disorder of objects and they don't need to
confuse themselves with this issue
 \end{quote}
 
\paragraph{
The above quote presents an interesting modification of Kant's original view. Following Kant Cassirer takes the issue of objethood seriously and continues to analyze it in Kantian (rather than Fregean) terms. However unlike Kant Cassirer treats the objecthood is an external issue with respect to the pure mathematics saying that the pure mathematics just like formal logic deals only with the ``order of concepts''. Thus he by and large accepts Russell's view on mathematics and displaces the problem of obecthood elswhere, namely into the sphere of application of mathematics to the sensual experience (in particular, within physics). Unlike Russell Cassirer doesn't try to reform mathematics but rather aims at describing it as a fact (``fact of science'' in Hermann Cohen's sense). So what Russell says about mathematics in 1903 Cassirer takes to be an adequate description of the current state of affairs in mathematics rather than a project aiming at reforming this discipline.
} 
    
\footnote{
Russell's logicism about mathematics is the claim according to which mathematics is wholly reducible to logic - in the sense that mathematics needs no other first principles except logical principles. A milder form of logicism about mathematics consists of the claim according to which first principles of mathematics contain logical principles plus some additional non-logical axioms like axioms of ZF. Historically such a mild logicism can be traced back at least to Frege's view about geometry \cite{Frege:1971}
Some people may argue that such a view does not qualify as logicism at all. Nevertheless I shall use in this paper the term ``logicism'' in such a wider sense that covers both Russell's strong logicism and the mild logicism just described. What Cassirer says about mathematics in the above quote in the main text does not commit him to the strong Russell's logicism but it does commit him to the mild one. 
} 

\section {Formal and Informal Bourbaki}

\paragraph{
It must be stressed that during the first decade of the 20th century things in mathematics and mathematical physics were changing so rapidly that no established pattern of doing mathematics and physics really existed at that time; it was a period when older patterns were competing with new nascent patterns. Although today we are still facing rapid conceptual changes the passed 20th century provides us with a historical distance allowing one to revaluate critically what Russell, Cassirer and other people who wrote about mathematics in 1900-ies. As seen from the beginning of the 21st century the situation appears to be the following. Russell's individual intellectual development from the (neo)Kantian critical philosophy to the (neo)Hegelean dialectical logicism to a Leibniz-style metaphysical  logicism supported by the new symbolic logic \cite{Russell:1946} represents a significant trend in the philosophy of the 20th century, which became particularly influential within the part of Anglo-Saxon philosophy that describes itself as Analytic. A mild form of logicism (see the last footnote) became a dominant trend in the philosophy of mathematics in spite of the fact many (if not most of) working mathematicians beginning with Poincar\'e disagree with it. Most mathematicians resist it by saying that they are simply not interested in philosophical and foundational issues in their professional lives.  Studies of formal axiomatic theories developed into an established field of mathematics but philosophical (in particular, metaphysical) aspects of these studies are seen by the most of the mathematical community as an alien element of the mathematical research. Like in Kant's time - and arguably even more than in Kant's time - today's mathematical community remains sensible to the difference between a mathematical argument and a philosophical speculation.   
}   

\paragraph{ 
Looking at patterns of today's mathematical reasoning one can not only recover the Euclidean scheme but also reaffirm Kant's argument according to which the objective character of mathematics depends on the way, in which mathematics constructs its objects.
}

\paragraph{ 
I shall demonstrate the relevance of this Kantian view in two steps: first, by analyzing recent informal mathematics and then by analyzing an example of wholly formalized theory. In both cases I shall use the example of Bourbaki's volume on Set theory (the first volume of Bourbaki's \emph{Elements} but consider two different editions of this work: for the formalized part I shall use the usual English edition \cite{Bourbaki:1968}, and for the informal part I shall use a draft version of this volume found in the Bourbaki's archive. While the final published version of the text contains a purely formal treatment of Set theory and aims at showing the possibility of a similar treatment of any other mathematical theory the draft presents Set theory in an informal manner, which quite accurately reflects the way in which this theory is used in the common mathematical practice (except the Set theory itself, as explained the previous Section). Let us begin with the draft.  
}

\footnote{\emph{El\'ements de la th\'eorie des ensembles}, R\'edaction 050, COTE BKI 01-1.4, available at mathdoc.emath.fr/archives-bourbaki/PDF/050\_iecnr\_059.pdf
}

After a philosophical introduction that explains the notion of mathematical theory the author introduces (informally) the concept of \emph{fundamental set} and the relation of membership between sets and their elements. (The author calls a set \emph{fundamental} in order to distinguish a well-defined set concept from a more general notion of collection.) Then the author introduces (with the usual informal notation) the concept of \emph{subset}, which is the subject to the following axiom:  

 \begin{quote}
Any predicate of type $A$ defines a subset of $A$; any subset of $A$ can be defined through a predicate of type $A$. 
 \end{quote}

(Predicate of type $A$ is a predicate $P$ such that for every element $a$ of set $A$ $P(a)$ has a definite truth-value. The subset $S$ of set $A$ defined by $P$ consists of such and only of such $a$ for which $P(a)$ is true.)  

Next Bourbaki introduces the concept of \emph{complement} of a given subset, of \emph{powerset} $P(A)$ of a given set $A$ (i.e., the set of all subsets of $A$); of \emph{union}, \emph{intersection} and \emph{cartesian product} of sets (described as \emph{operations} on sets), of \emph{relation} and \emph{function} between sets. Having these basic concepts in his disposal the author says:

\begin{quote}
In any mathematical theory one begins with a number of fundamental sets, each of which consists of elements of a certain type that needs to be considered. Then on the basis of types that are already known one introduces new types of elements (for example, the subsets of a set of elements, pairs of elements) and for each of those new types of elements one introduces sets of elements of those types. 
\\  
So one forms a family of sets constructed from the fundamental sets. Those constructions are the following: \\
1) given set $A$, which is already constructed, take the set $P(A)$ of the subsets of $A$;\\
2) given sets $A$, $B$, which are already constructed, take the cartesian product $AxB$ of these sets.\\
The sets of objects, which are constructed in this way, are introduced into a theory step by step when it is needed. Each proof involves only a finite number of sets. We call such sets \emph{types} of the given theory; their infinite hierarchy constitutes a \emph{scale of types}.   
\end{quote}
  
On this basis the author describes the concept of  \emph{structure} in the following way:  
\begin{quote}
We begin with a number of fundamental sets: $A, B, C, ... , L$ that we call \emph{base sets}. To be given a \emph{structure} on this base amounts to this: \\
1) be given properties of elements of these sets;
2) be given relations between elements of these sets;
3) be given a number of types making part of the scale of types constructed on this base;
4) be given relations between elements of certain types constructed on this base;
5) assume as true a number of mutually consistent propositions about these properties and these relations. 
\end{quote}

\paragraph{
What I want to stress is the fact that principles of building mathematical theory described in the Bourbaki's draft are not so different from Euclid's: Bourbaki like Euclid begins with principles of building mathematical objects but not with axioms. Axioms (in the modern sense of the term) appear only in the very end of the above list (the 5th item). While for Euclid the basic data is a finite family of \emph{points}  and the rest of the geometrical universe is constructed from these points by Postulates for Bourbaki the basic data is a finite family of \emph{sets} and the rest is constructed as just described. While for Euclid the basic type of geometrical object is a \emph{figure} for Bourbaki the basic type of mathematical object is a \emph{structure}. In both cases the constructed objects come with certain propositions that can be asserted about these objects without proofs because they immediately follow from corresponding definitions. In both cases the construction of objects is a subject of certain rules but not the matter of a mere stipulation. An infinite cyclic group (to take a typical example) construed as a \emph{structure} along the above lines qualifies as an \emph{object} just like a regular triangle in Euclid.
}

\footnote{In order to continue this analogy one may compare the notion of \emph{isomorphism} of structures in Bourbaki with Euclid's notion of \emph{congruence} of figures.}

Remind however that we are talking so far about a draft but not about the finial published version of the Bourbaki's volume on Set theory. The published formalized version of the above is drastically different. Here Bourbaki proceeds as follows. The first chapter of the treatise, which has the title \emph{Description of Formal Mathematics}, begins with an account of \emph{signs} and \emph{assemblies} (strings) of signs provided with a definition of mathematical theory according to which a theory

\begin{quote}
... contains rules which allow us to assert that certain assemblies of signs are \emph{terms} or \emph{relations} of the theory, and other rules which allow us to assert that certain assemblies are \emph{theorems} of the theory.  
\end{quote}

\paragraph{
Then follows a description of operations that allow for constructing new assemblies of signs from some given assemblies; the simplest operation of this sort is the \emph{concatenation} of two given assemblies $A, B$ into a new assembly $AB$.  On such a purely syntactic basis Bourbaki introduces some logical concepts necessary for a formal axiomatic treatment of Set theory. Although  Bourbaki's version of axiomatic Set theory is not identical to ZF it is similar in its character; the differences between the two ways of formalizing Set theory are not relevant to the present discussion and I leave them aside. Instead, I shall stress differences and similarities between the informal Bourbaki's draft on Set theory and the published volume, where Set theory is treated formally.
}
\paragraph{
In the published formalized theory (unlike the theory found in the unpublished draft) every basic set-theoretic \emph{construction} is represented as an assembly of signs (string of symbols) which, informally speaking,  states (as an axiom) or proves that given certain sets there exist some other sets (like in the case of \emph{pairing}: given sets $x, y$ there exist another set $z = \{x, y\}$ that has $x$ and $y$ as its elements). Thus in the formal theory sets are no longer  \emph{constructed} in anything like Kant's sense ; instead one assumes the \emph{existence} of certain sets and using this assumption proves the existence of certain other sets.   
}

This is, of course, a very general feature of the modern axiomatic method, which I have already stressed above talking about modern formal representations of the Euclid-style geometrical reasoning. Hilbert and Bernays  \cite{Hilbert&Bernays:1934} describe this difference between their own approach and Euclid's approach in the following words:  

\begin{quote}
Euclid does not presuppose that points or lines constitute any fixed
domain of individuals. Therefore, he does not state any existence
axioms either, but only construction postulates.
\end{quote}

\paragraph{
However as a matter of fact Bourbaki's formal Set theory just like its informal prototype begins with certain basic constructions but not with axioms! In that respect there is no essential difference between the formal and the informal approaches, and the Kantian view on mathematical objecthood applies equally well in both cases. The difference is more specific and concerns the \emph{types} of involved constructions.  While in the informal draft the relevant constructions are set-theoretic (subset, complement, powerset, union, intersection, cartesian product, etc.) in the published formalized version of the theory the relevant constructions are  \emph{syntactic} (assembly of signs, concatenation of assemblies, substitution, etc.). In order to understand this difference we need to discuss first the role of symbolism in mathematics. 
}

\section {The Role of Symbolism}

The constructive aspect of formalization just stressed remains unrevealed in Russell's logicist philosophy of mathematics, which uses the symbolic logic as a magic tool that helps to realize the Leibnizian dream but doesn't explain what does it mean for logic to be \emph{symbolic}. Hilbert and Bernays in the introductory part of their \cite{Hilbert&Bernays:1934} provide the following explanation of this issue. First they remark that for developing the elementary arithmetics    

\begin{quote}
what we essentially need  is only that the numeral 1 and the suffix 1
are intuitive objects which can be recognized unambiguously
\end{quote}

(meaning that every natural number can be represented as a string of 1s.) Then they

\begin{quote}
also briefly characterize the elementary conceptual viewpoint in algebra. ... The objects of the theory are again certain figures, namely the polynomials constructed with the help of the symbols $+, - , \bullet $ and parentheses.
\end{quote}

Finally, considering the formalization of Calculus they say

\begin{quote}
When the usual Calculus is formalized  (i.e. when its presuppositions and inferences are translated into initial formulas and rules of deduction), then a proof in Calculus becomes a  succession of intuitively comprehensible processes. ... Then, in principle, we have the same situation as in our treatment of the elementary arithmetic.
\end{quote}

\paragraph{
Hilbert and Bernays talk here about mathematical \emph{objects} and about their intuitive character in the full accordance with Kant. The purpose of formalization, in their understanding, is \emph{not} getting rid of objecthood and intuition in mathematics but rather a replacement of the problematic objecthood and the problematic intuitions of modern mathematics (which involves non-Euclidean spaces, infinite sets and the like) by the apparently non-problematic objecthood and intuition related to manipulations with finite strings of symbols according to some precisely specified rules. 
}

\paragraph{
Let me now make another historical digression and remark that mathematics began to use symbolic means at a very early stage of its historical development. The history of mathematical notation \cite{Cajori:1929} provides, in my view, an appropriate context for considering the modern formal mathematics. For this purpose I shall describe here briefly two historical examples and then use them for a further analysis of the Hilbert-style formalization.   
}

\paragraph{
My first example is Euclid's mathematical notation.  As we have seen Euclid uses written means of three different types: (i) the usual phonetic alphabetic writing by which he communicates his Greek prose, (ii) diagrams and (iii) a special alphabetical mathematical notation (like in the case of $ABC$ used as a name of triangle). This notation links (i) with (ii) in the following way: a name like $ABC$  can be inserted into the Greek prose and can be equally ``read of'' the corresponding diagram; in this way this notation allows for referring to the diagram in the prose. This notation remains today in use in the school geometry; although Euclid uses it systematically both in geometrical and in arithmetical Books of his  \emph{Elements} in the following discussion I leave here Euclid's arithmetic aside.  
}

\footnote{(i) Although the early MSS of Euclid's \emph{Elements} contain no diagram there is good reason to believe that Euclid and other Greek mathematicians did use diagrams in their practice, most probably by drawing them on a specially prepared sand \cite{Netz:1999}\\
(ii) It is sometimes said that in the arithmetical Books of the  \emph{Elements} Euclid represents numbers by geometrical straight lines. This, in my view, is not accurate. Indeed Euclid represents straight lines and numbers by diagrams of the same type and he uses in both cases the same notation. It doesn't follow that he represents numbers by straight lines. Any mathematical diagram and any mathematical notation requires a convention \emph{what} the given diagram and the given notation represent. The same drawings and the same symbolic expression may represent in mathematics very different things.
}

\paragraph{
Euclid's geometrical notation has an implicit syntax that makes this notation appropriate for its purpose. In particular, it fits Euclid's rules about constructing geometrical objects. Constructing the name ``$AB$'' from the names ``$A$'' and ``$B$''  mimics drawing the (segment of) straight line $AB$ between the given points $A$ and $B$. The convention according to which names ``$AB$'' and ``$BA$'' always refer to the same straight line reflects the assumption according to which drawing a straight line from $A$ to $B$ and drawing a straight line from $B$ to $A$ always brings the same result. This syntax extends to names for polygons as follows. Every two-letter segment of the name of a given polygon is the name of a side of this polygon; by concatenating the first and the last letter of the name of the polygon one also gets the name of a side of this polygon.  This convention and the assumption according to which a polygon is uniquely determined by its sides imply that a cyclic permutation of the letters in a name of polygon remains the referent of this name invariant. For example if word ``$ABC$'' names a triangle than word ``$CBA$'' names the same triangle. Another important correspondence between syntactic rules and constructive rules in Euclid concerns the \emph{sum} and the \emph{difference} of numbers and magnitudes that is implicitly defined by Axioms 2 and 3 of the  \emph{Elements}. In particular, these operations make sense for straight lines (albeit they are not defined for \emph{any} pair of straight lines). Consider two straight lines $AB$ and $BC$ such that point $B$ lies on $AC$. In this case the sum of $AB$ and $BC$ is well-defined and $AB + BC = AC$. Summing up of $AB$ and $BC$ is, of course, a geometrical rather than syntactic operation: it takes two geometrical objects and produces out of them a new object. But it is mimicked on the syntactic level through the following obvious rule: given two names $UV$ and $XY$ such that $V = X$ ``contract'' them and get new name $UY$. Although such a syntactic contraction does not fully reflects the corresponding geometrical operation (in particular, it doesn't allow one to see whether or not the corresponding geometrical operation is well-defined) it makes the syntax convenient for dealing with this geometrical operation.    
}

\paragraph{
The next historical example is  suggested by Hilbert and Bernays' reference to the elementary algebra in the above quote. In the 17th century Descartes and his followers established a new symbolic algebraic notation, which is still in use (in only slightly modified form) in the elementary school algebra and which is in many respects similar to the notation used in the modern abstract algebra. The new notation had a tremendous effect on mathematics recently described by Serfati \cite{Serfati:2005} a ``symbolic revolution''. The symbolic calculus known today as the elementary algebra in the 17th century was understood by a number of influential writers as a general science of magnitude; the relevant notion of magnitude generalized upon natural numbers and geometrical magnitudes known from Euclid and other traditional sources. Consider the following interesting passage from Arnauld \cite{Arnauld:1683}:
}
\begin{quote}
 What cannot be multiplied by its nature can be multiplied through a mental fiction where the truth presents itself as certainly as in a real multiplication. In order to learn the distance covered during 10 hours by one who covers 24 lieu per 8 hours I multiply through a mental fiction 10 hours by 24 lieu that gives me the imaginary product of 240 hour times lieu, which I divide then by 8 hours and get 30 lieu.  By the same mental fiction one multiplies a surface by another surface  even if the product has 4 dimensions and cannot exist in nature. One may discover many truths through multiplications of this sort. \\ 
 People say that this is because the imaginary products can be reduced to lines. ... But there is no evidence that [relevant] proofs  depend on those lines, which are in fact wholly alien to them.  (p. 38-39) 
\end{quote}

\paragraph{
 Traditionally the product of two straight lines is construed as a rectangle having these given lines as its sides; the product of three lines is a solid but in order to form products of four and more linear factors in this geometrical way one needs higher dimensions, which according to Arnauld ``cannot exist in nature''. Nevertheless he is ready to consider such higher-dimensional geometrical constructions as useful fictions on equal footing with products of distances by time intervals, which don't have any immediate physical interpretation either but are demonstratively useful for calculations. In the last quoted paragraph Arnauld refers to Descartes' proposal to construct the multiplication of geometrical magnitudes differently, so the product of straight lines is again a straight line; in modern words Descartes' definition of multiplication of straight lines makes this operation algebraically closed. (Descartes uses an auxiliary line 1 as a unit and then considers similar triangles with sides $1, a, b, c$ such that $\frac{1}{a} = \frac{b}{c}$, which gives him the wanted definition of $c = ab$, see \cite{Descartes:1886}) Arnauld finds this trick artificial and unnecessary. What makes him confident about higher-dimensional geometrical products and quasi-physical units like hour times lieu is the symbolic algebraic calculus, to which these problematics notions are associated. Within this calculus the product of four factors $p = abcd$ behaves like a product of two or three factors. One cannot easily imagine the product of two surfaces (just like one cannot give a physical sense to hours times lieu) but one can easily concatenate the string $ab$ and the string $cd$ and think of this operation as multiplication of two surfaces. Descartes' alternative definition of the geometrical product aims at providing a ``clear and distinct''  intuitive underpinning of this operation that avoids the talk of higher dimensions. Arnauld finds Descartes' construction of product unnecessary because in his eyes the symbolic calculus provides such an intuitive underpinning by itself. The ``proofs'' that Arnauld mentions in this context are nothing but symbolic calculations.  
}

It is instructive to compare Arnauld's treatment of algebra with the following Kant's remark about this subject:

\begin{quote}
[M]athematics does not merely construct magnitudes (quanta), as in geometry but also mere magnitudes (quantitatem), as in algebra, where it entirely abstracts from the constitution of the object that is to be thought in accordance with such a concept of magnitude. In this case it chooses a certain notation for all construction of magnitudes in general ....  and thereby achieves by a symbolic construction equally well what geometry does by an ostensive or geometrical construction (of objects themselves), which discursive cognition could never achieve by means of mere concepts. \cite{Kant:1999}, A717
\end{quote}

\paragraph{
In this quote Kant presents a view on the elementary symbolic algebra, which can be called \emph{formal} (in the modern rather than in Kant's proper sense of the term). According to this view ``true'' objects of symbolic algebra are symbolic constructions, which represent ``mere magnitudes'' (or ``pure magnitudes'') just like constructed geometrical triangles represent the general concept of triangle. From the intuitive point of view such symbolic constructions are at least as clear and distinct as traditional geometrical constructions and, arguably, are even more clear and more distinct. In a sense these symbolic constructions are even ``more concrete'' than traditional geometrical constructions: as Hilbert and Bernays tell us in the above quote such symbolic construction are as much concrete as a representation of natural number $n$ by $n$ vertical bars. At the same time  mathematical \emph{concepts}  represented by these symbolic constructions are more general and more abstract than traditional geometrical concepts: those are concepts of \emph{abstract magnitudes} rather than concepts of circle, triangle and similar familiar things.  
}

\paragraph{
Such a formal understanding of algebra is apparently not alien to Arnauld either. As far as algebraic calculations are concerned he readily relies on symbolic procedures and finds it unnecessary to justify these procedures by any geometrical considerations. Yet Arnauld's words (as well as commonly known facts of the history of mathematics of the last three centuries) make it clear that the purely formal view on the symbolic algebra does not fully express the significance of this mathematical discipline.  Another part of the truth is that the symbolic algebraic constructions allow for multiple geometrical and arithmetical interpretations and moreover may motivate the development of new geometrical and arithmetic concepts like that of higher-dimensional space or that of complex number. This feature of the symbolic algebra of the 17th century allowed Descartes and his follows to apply the new algebra to geometry, which gave rise to what we call today \emph{analytic geometry}. The dialectical interplay between the abstract algebra, on the one hand, and geometry and arithmetic, on the other hand, still remains one of the most important driving forces of the mathematical progress  (think of recent achievements in the \emph{algebraic geometry}). This shows that the \emph{contentual} aspect of symbolic algebra has been at least as much important throughout the long history of this discipline as its formal aspect.  
}

\paragraph{
Both Euclid's geometrical notation and the symbolic algebraic notation represent mathematical operations and mathematical constructions by certain operations with symbols and certain symbolic constructions. The algebraic notation does this, of course, in a more systematic way, which makes possible a purely formal reading of this notation. (The proto-syntax of Euclid's notation makes this notation convenient but doesn't allow one to use this notation \emph{without} diagrams except very simple cases.) Let us now look more precisely at the analogy between the elementary symbolic algebra and symbolic formal calculi suggested by Hilbert and Bernays in the above quote.  At the first glance the analogy seems to be precise. In both cases we have symbolic constructions provided with precise syntactic rules, on the one hand, and some problematic reasoning involving confusing intuitions (higher geometrical dimensions, imaginary magnitudes, infinite sets, etc.), on the other hand. In both cases the symbolic constructions and the associated syntax is used for clarifying the problematic reasoning. In both cases a \emph{formal} approach is appropriate. No matter how one imagines 4-dimensional objects, quasi-physical units like hours times lieu and infinite sets, one knows how to manipulate with these things using symbols and precise rules of manipulation with these symbols. So as long as 4-dimensional objects, quasi-physical units and infinite sets are represented by symbols in an appropriate symbolic calculus one is on a secure ground (provided, of course, that the calculus itself is secure).         
}

\footnote{Compare what MacLaurin says about the usefulness of algebra in his \emph{theory of fluxions}:  
\begin{quote}
The improvement that have been made by it [the doctrine of fluxions] ... are in a great measure owing to a facility, conciseness, and great extend of the method of computation, or algebraic part. It is for the sake of these advantages that so many symbols are emplyed in algebra.  ... It [algebra] may have been employed to cover, under a complication of symbols,  obstruse doctrines, that could not bear the light so well in a plain geometrical form; but, without doubt, obscurity may be avoided in this art as well as in geometry, by defining clearly the import and use of the symbols, and proceeding with care afterwards. \\
(\emph{A Treatise of fluxions}, quoted by \cite{Cajori:1929},  v2, p. 330) 
\end{quote}
}

\paragraph{
As it is well known the symbolic calculi appropriate for the purpose turned to be, generally, less secure than Hilbert and his followers hoped for. I would like, however, to stress a very different aspect of the problem by pointing on a \emph{disanalogy} between the symbolic algebra, on the one hand, and the symbolic calculi used for formalization of mathematics, on the other hand. As we have seen the symbolic algebra represents certain problematic geometrical constructions like the product of two surfaces by simple symbolic constructions like the concatenation of strings $ab$ and $cd$. The \emph{informal} set-theoretic notation like one used in the Bourbaki's draft works similarly: it helps to represent set-theoretic constructions with certain symbolic constructions as in the case when the set of four elements $a, b, c, d$ is denoted as $\{a, b, c, d\}$. However the properly \emph{formal} syntax like one used in the Bourbaki's published volume on Set theory or in ZF works quite differently. In this case symbolic constructions represent not set-theoretic constructions but \emph{inferences} of some propositions from some other propositions. As we have already seen the Hilbert-style formal approach to Set theory assumes that in mathematics (albeit not in the \emph{meta-mathematics}) the talk of \emph{construction} can be safely replaced by the talk of \emph{existence}. Thus within a Hilbert-style formal logical framework syntactic manipulations with symbols no longer represent manipulations with mathematical objects as this is done in the symbolic algebra (where manipulations with symbols represent manipulations with geometrical and arithmetical objects or with pure magnitudes) but rather represent manipulations with propositions (which under the intended interpretation tell us something \emph{about} certain mathematical objects). In particular, syntactic operations described in the beginning of the Bourbaki's volume on Set theory do not represent set-theoretic operations informally described in the Bourbaki's draft.  
}

\paragraph{
I claim that this feature of the Hilbert-style formal axiomatic method makes this method inappropriate for doing mathematics, and I hold it responsible for the enduring gap between the formal and the informal mathematics exemplified by the two Burbaki's texts discussed above. This is, in my view, a reason why the formalization of mathematics so far did not have an effect comparable with the symbolic revolution of the 17th century but isolated itself within a specific field of research called by its proponents ``foundations of mathematics'' (in spite of the fact that most of working mathematicians pay little interest to this research and don't use its results in their work). Examples from today's mainstream ``informal'' mathematics strongly suggest that mathematics did \emph{not} change its nature from Euclid's times so radically that the object-building eventually became a non-issue and was replaced by a purely speculative stipulation of abstract entities with some desired properties relations. Today's mathematics just like Greek mathematics constructs its \emph{objects} like groups, topological spaces, toposes and whatnot. Just like Greek mathematics today's mathematics constructs its objects following certain \emph{rules}, some of which apply across the whole of mathematics and some of which are specific for a given mathematical discipline. Just like Greek mathematics today's mathematics involves \emph{doing and showing} but not just showing and proving. True, informal constructive rules used in the contemporary mathematics in many cases call for more precise specifications. I can see no general reason why such more precise specifications cannot be given but must be replaced by appropriate existential propositions. 
}

\paragraph{
If we now look more precisely where exactly lies the problem we can see that it has to do with the common \emph{informal} notion of axiomatic theory discussed in the First Part of this paper but doesn't concern more specific and more technical features of the modern axiomatic method.   
}

\section{Mathematical Constructivism}
\paragraph{
In spite of the dominant position of (weakly) logicist approaches to mathematics during the most of the 20th century the Kantian view according to which mathematics should construct its object rather than just stipulate them, was never wholly forgotten and motivated a significant amount of mathematical and philosophical work that constitutes a trend in mathematics and its philosophy called \emph{constructivism}. This is a rather diverse movement that I shall not purport to describe here in detail; for a concise review of the mathematical constructivism of the late 19th and the 20th century (which however doesn't cover recent developments) and further references see \cite{Troelstra:1991}. Here I provide only few remarks about it, which are relevant to my purpose.    
}

\paragraph{
Rather unfortunately, in my view, the constructivist movement from the very beginning took a conservative bend and began a fight against then-new ways of mathematical thinking including the set-theoretic thinking. This tendency can be traced back to Kronecker who required every well-formed mathematical object to be constructible from natural numbers   (famously saying that only natural numbers are God's creation while everything else in mathematics is a human handicraft) and on this basis didn't admit infinite sets. More recently Bishop was inspired by similar ideas (and in particular by Kronecker's works). Brouwer's \emph{intuitionism} (which qualifies as a form of constructivism) also put rather severe restrictions on his contemporary mathematics as well as on essential parts of earlier established mathematical results. Even those constructivists who like Markov tried to develop constructive mathematics as a special part of mathematics rather than reform mathematics as a whole understood the notion of mathematical construction very restrictively and almost exclusively in computational terms.   
}

\paragraph{
Since in this paper I also push the Kantian line and stress the significance of constructions and objecthood in mathematics my present view also qualifies as constructivist. However the target of my constructivist critique is not some particular mathematical notion like that of infinite set (of the countable or a higher cardinality) that may look suspicious from a constructive point of view but the modern axiomatic method itself. I claim that a mathematical theory needs explicitly given \emph{constructive rules}, which cannot be, generally, reduced to rules of logical inference (whether the inference is treated formally or not). I make so far no judgement as to which sorts of constructions are admissible and which are not. As long as we are talking now about a constructive method of theory-building in general the discussion about particular ways of constructing mathematical objects is inappropriate just like  the discussion on particular axioms is inappropriate when one discusses generalities concerning the axiomatic method. The Kroneckerian belief that only arithmetical constructions are ``real'' is not very unlike the belief that only figures constructible by ruler and compass are real. Any constructive thinking in and about mathematics must definitely avoid such a dogmatic attitude just like Hilbert's axiomatic thinking avoids the dogmatic attitude towards particular axioms and focuses instead on questions concerning the relative consistency of axioms and the like.  
}

\paragraph{
The example of Euclid's  \emph{Elements} shows how a perfectly constructive theory can be built in a systematic manner. As we have seen this theory is not built as an axiomatic theory in the modern sense of the word. The main difference is that Euclid's theory includes such non-propositional constructive principles as \emph{postulates}. I claim that such non-propositional constructive principles are needed in modern mathematical theories too. }

\paragraph{
The idea that constructive postulates in a formal axiomatic can be always replaced by existential axioms may look innocent but it is not. Let's see again what is gained and what is lost with such a replacement (leaving now aside metaphysical issues about being and becoming). In such a propositional setting the problem of constructibility of object $O$ on the basis of a given set of postulates $P$ reduces to the problem of provability of certain existential proposition $E$ (that says that $O$ exists) from a set of assumptions $A$ (that includes axioms of the given theory and eventually some additional hypotheses). In other words one should provide proof $B$ that uses assumptions $A$ and brings $E$ as a conclusion. Since in a formal setting $B$ is itself a (syntactic) construction one thus reduces the question of constructibility of $O$ to the question of constructibility of $B$. The constructibility of $B$, of course, depends of logical rules of inference $R$, which are used in the given axiomatic setting. Beware that $R$ and $P$ are sets of rules but not sets of propositions.  
}
\paragraph{
Now if object $O$ and postulates $P$ are described only informally like in the Euclid's  \emph{Elements} or in the Bourbaki's draft mentioned above one may argue that the reformulation of the constructibility problem in terms of formal proof $B$ makes sense because the syntactic object $B$ is simpler than $O$ and that rules $R$ are more precise and better determined than rules $P$. ``Simpler'' and ``better determined'' basically means that $B$ and $R$ has a finitary combinatorial nature. However such a reduction of one constructibility question to another constructibility question (which is better posed and admits a more precise answer) is also achieved by a somewhat similar but still a different means, namely by the use of symbolic algebra. And at least as far as the question of constructibility by ruler and compass is concerned the use of symbolic algebra proves more effective than the use of symbolic logic because algebraic symbolic constructions mimic geometrical constructions in a sense, in which symbolic constructions of a formal counterpart of Euclid's geometry do not.      
}

\footnote
{ See Tarski's classics \cite{Tarski:1959} where the author considers several possible axiomatizations of Euclidean geometry and study their meta-theoretical properties. A relevant axiomatic theory, which takes into account Euclid's constructibility by ruler and compass turns to be both  incomplete (in the sense that some propositions formulated in the language of this theory are true in certain models of the theory but false in other models) and not decidable (there exist no general procedure for proving or disproving propositions formulated in the language of this theory). So at least in this particular case the reduction of informal mathematical constructions to formal symbolic constructions doesn't help one to answer questions about constructibility of particular objects in the original theory.
}

\paragraph{
The principle reason why one may be not satisfied with the algebraic solution lies beyond the issue of constructibility as such. One may argue that even if the informal symbolic algebra indeed clarifies the issue of geometrical constructibility it leaves imprecise the general logical structure of the argument. A formal axiomatic approach aims at reducing the whole argument - and not just the part of the argument that involves object-construction - to some sort of symbolic calculation. This is why the Kantian distinctions between the logic of concepts (i.e., formal logic in Kant's sense) and the logic of intuitive representations of these concepts or otherwise the logic of corresponding \emph{objects} (transcendental logic) does not make much sense in the context of formal theories. On the one hand, a formal proof (in the modern sense of ``formal'') can be seen as an argument based on certain assumptions and applying certain logical rules, i.e., as a formal proof in Kant's sense. One the other hand, it can be seen as an intuitively transparent  symbolic construction, which is a subject to precise constructing rules that can be described in terms of transcendental logic. Strictly speaking none of these two readings is relevant because in the given context the very difference between construction and calculation, on the one hand, and purely conceptual non-objectual reasoning on the other hand, wholly disappears. This seems me to be the core of the neo-Leibnizian  project in foundations of mathematics of the 20th century, which purports to replace the Kantian critical approach. Its description along the above lines is found in \cite{Carnap:1931}.
}

\paragraph{
This idea may look indeed very modern (in spite of its quasi-traditional metaphysical aspect that I stressed in the previous section) and very radical (against the Kantian background). There are however two reasons why I believe that it is fundamentally flawed and and must be rather abandoned. 
}
\paragraph{
The first reason is practical. A century is a sufficient historical laps for judging whether a particular project of reforming mathematics has been successful or not. An analysis of the current mathematical practice suggests that the project of reforming mathematics along the neo-Leibnizian lines has not be successful. Working mathematicians still distinguish between objects that they construct and mathematical facts about these objects that they formulate and prove. Noticeably this concerns both the mainstream informal mathematics and studies in formal logic and related areas. In the informal mathematics people widely apply various forms of constructions and calculations but they also \emph{show} how and why these constructions and these calculations are relevant to their arguments. In the formal mathematics people do the same but focus their research on the specific case of symbolic constructions and specific symbolic calculations called \emph{logical}. In the formal mathematics just like in the informal one people not only calculate but also reason about what they calculate.  Thus neither the formal nor the informal mathematics provides an evidence of the neo-Leibnizian fusion of reasoning with calculation. True, mathematical proofs in some cases reduce to calculations. However such a reduction always makes part of a given proof. So this fact does not count as an evidence of the neo-Lebnizian fusion. Although there is always a space for arguing that the current mathematical practice is ill-founded and needs to be reformed I believe that the feedback from practice to theory is essential. This is why I take the above practical argument seriously and am not ready to devaluate it as ``merely pragmatic''.           
}

\paragraph{
The second reason is theoretical and it lies beyond the pure mathematics itself. It needs a special discussion. 
}

\section{Galilean Science and ``Unreasonable Effectiveness of Mathematics''}
\paragraph{
It is hardly controversial that mathematics deals with forms of possible human experience;  in its simplest and most general form this claim is simply tantamount to saying that mathematics applies across a wide range of human practices. Today this is even more true than it was in Kant's time: crucial technologies, on which depend our well-being, in many ways depend on mathematical considerations and cannot be sustained and further developed without mathematical expertise; mathematics today makes part of any engineering education. In Kant's time the only properly mathematized science was (Newtonian) mechanics; the following progress of science in the 19th century has brought us to the point when every physical theory deserving the name has a mathematical aspect. Today physics and chemistry are mathematized and the mathematization of biology is in progress. Using mathematical models also becomes an usual practice in social sciences. 
}

Let me now be more specific and ask which general forms of experience are relevant to today's science and technology. This question is obviously yet too large and too general to be answered in any detail here. I would like however to stress only the following point. At least since Galileo's times science practices an active intervention of humans into the nature through experiments rather than a passive observation and description of the observed phenomena. So we are talking now about the mathematically-laden science where mathematics serves for guiding human interactions with the environment rather than simply for describing how this environment appears to our senses. As van Fraassen \cite{Fraassen:1980} puts this

\begin{quote}
The real importance of theory, to the working scientist, is that it is a factor in experimental design. (p. 73)
\end{quote}

Thus mathematical forms of possible experience relevant to the modern Galilean science are forms of such possible \emph{interactions} with the environment rather than only linguistic and logical forms that allow for spelling out some plausible hypotheses about the world and deriving from them some consequences according to certain rules. The forms of the latter kind may be sufficient for developing a speculative science along the older scholastic pattern but they are certainly not sufficient for developing the modern mathematically-laden science and the modern mathematically-laden technology.           
 
 \paragraph{ 
As far as the pure mathematics is conceived as a domain of abstract logical possibilities the fact that mathematics proves ``unreasonably effective'' \cite{Wigner:1960} in its applications to empirical sciences and technology remains a complete mystery. The mystery is dissolved as soon as one observes that mathematics explores not everything that can possibly \emph{be} the case (which is a hardly observable domain unless one delimit the sense of ``possibly'' in one way or another) but rather what we can possibly \emph{do} within the limits of our human capacities (which are steadily growing with the progress of science and technology). What are these limits is a tricky question. On the one hand, mathematics systematically ignores certain apparent limits by exaggerating relevant capacities: this is usually called the \emph{mathematical idealization}. For example, mathematicians pretend that they can count up to $10^{10^10}$ just as easily as up to 10 or that they can draw a straight line between two stars just as easily as they can draw a line between two points marked on a sheet of paper. This strategy usually works until the point where the empirical constraints become pressing and people invent new mathematics that takes these constraints into account as this, for example, happened when people realized that the old good Euclidean geometry is not appropriate for describing the physical space at large astronomical scales (in spite of the fact that it still works amazingly well at the scale of a planetary system like ours!). One the other hand, it also happens that in a real experiment people observe what in terms of the assumed mathematical description of the given experiment qualifies as impossible as this happened in the Michelson-Morley experiment supposed to measure parameters of the ether flow around the Earth. In such cases people say that the assumed mathematical description (and hence the corresponding physical theory) is wrong and look for a new one. Sometimes the suitable mathematics can be found in a nearly ready-made form and only used for building a new physical theory but sometimes in order to fix the problem one needs to develop the appropriate mathematics from the outset as this happened in the history of the electro-magnetism, for example. This picture suggests the view on the pure mathematics as a proper part of the Galilean science. However complicated the dialectic of interaction between the purely mathematical part and the pure empirical part of the Galilean science may be (in fact these parts hardly exist in pure states) there is nothing unreasonable in it. On the contrary, this dialectics shows us how exactly the Galilean scientific reason works.    
} 

\paragraph{
Russell's neo-Leibnizian logicism about mathematics promises nothing more and nothing less than that: to make mathematics a part of logic, so that any mathematical form of possible experience turns into the form of a \emph{proposition} (and forms of logical inference of propositions from some other propositions), which may eventually refer to some experience. Russell's view is quite radical in this respect, and many people including Hilbert who were directly involved into reforming mathematics on the basis of new logic in the beginning of the 20th century didn't share Russell's philosophical views. Anyway, as I have already argued, a weaker form of the neo-Leibnizian approach (that I call the \emph{weak} logicism in order to distiguish it from Russell's radical logicism) is intrinsic to the formal axiomatic method in Hilbert's sense of the term (which nowadays has became standard). Even if forms of possible experience delivered by a formal axiomatic theory do not qualify as \emph{logical} forms in the precise sense of the term they are still can be described as forms of possible empirical propositions rather than forms of empirical interactions or anything else.         
}

\paragraph{
As far as we want to continue to develop the Galilean science (and the technology connected to this type of science) our mathematics must provide for it forms of possible empirical interaction rather than just forms of propositions. In other words it must provide forms appropriate \emph{doing} various things in the world but not only forms for talking about this world. Since formalized mathematical theories are not appropriate for this job we need to learn how to build mathematical theories differently.}

\footnote{In certain situations one is in a position of doing things by saying some other things, for example, by giving orders to other people. Physical theories allow for a similar control over natural phenomena, which is the basis of the modern science-based technology. However just like in the case of human society such a control cannot be established only by words and requires controlling mechanisms of different sorts. The 
idea that such a control has been already established for us in advance by an external agency, so that the mysteries of the world disclose themselves as soon as one finds right words for talking about them, belongs to magic rather than science. 
}

\paragraph{
As Kant shows in great detail the traditional geometry and his contemporary algebra are useful in the Galilean science because these mathematical theories are \emph{constructive} in the sense that they involve rules for constructing their objects (explicitly or implicitly). Today we can hardly hope, of course, to get a new mathematical theory that would allow for identifying a physical object with a mathematical object in the same way, in which one may identify (modulo the mathematical idealization), say, a planet with an Euclidean sphere. Today people doing particle physics describe particles using the mathematical group theory and manipulating with particles in experiments using a special hi-tech equipment;  they don't expect that mathematical manipulations with groups would map their experimental manipulations in a direct way. Nevertheless  the constructive character of the mainstream informal mathematical practice, which I have stressed earlier in this paper, still helps physicists and other scientists to design their experiments and their equipments.  Scientists make up mathematical models of their experimental systems identifying (modulo the usual mathematical  idealization) their experimental systems with the models of these systems, then manipulate both with the models (theoretically) and with the experimental systems (in real experiments) and see whether the manipulations of both sorts work coherently. This is, of course, an oversimplified picture of the scientific experiment (for more details see \cite{Fraassen:1980}) but it is sufficient for seeing that the possibility to establish a correlation between mathematical manipulations, on the one hand, and experimental manipulations, on the other hand, remains essential for today's mathematically-laden experimental science. 
}
\paragraph{
Such a correlation cannot be possibly established when the only type of mathematical objects available for manipulation are \emph{syntactic} objects. I don't want to say that manipulations with syntactic objects cannot be useful in physical experiments  - in fact they can - but claim that the gap between syntactic manipulations and physical experimental manipulations needs to be filled by manipulations with mathematical objects (like groups and whatnot). The Hilbert-style formalized mathematics cannot provide this because it disqualifies the notion of manipulation with mathematical objects as naive and ill-founded and allows only for manipulations with syntactic symbolic constructions. Within such a formal approach one is allowed to say and proof various things about mathematical objects but not allowed to touch these objects and moreover to manipulate with them. Even if it is possible to develop some sort of sterile mathematics within such a restrictive setting such mathematics would be useless in the modern empirical science and moreover in the modern technology.   
}

\paragraph{
This shows that the constructive aspect of mathematics (in the relevant sense of ``constructive'' explained above) is indeed essential it must be taken taken into account by any useful method of theory-building applicable in mathematics.  Mathematics certainly needs formal languages and logical rules applicable to propositions expressed in these languages but it also needs rules of a different sort, namely rules that regulate constructions of and operations with non-syntactical mathematical objects. Just like logical rules such constructive non-logical rules can be described syntactically as rules about transformation of some symbolic constructions into some other symbolic constructions. So what I have in mind is a development of the modern axiomatic method rather than its complete replacement by something else. I believe that the example of Euclid's  \emph{Elements}, which inspired Hilbert about a century ago, can be once again useful for this purpose.
} 

\bibliographystyle{plain} 
\bibliography{doshow}

\end{document}